\documentclass{amsart}
\usepackage{amsfonts}

\usepackage{graphicx}

\input{tcilatex}
\input tcilatex

\begin{document}
\title[Pareto genealogies in a Poisson evolution model with selection]{Pareto genealogies arising from a Poisson branching evolution model with
selection}
\author{Thierry E. Huillet}
\address{Laboratoire de Physique Th\'{e}orique et Mod\'{e}lisation\\
CNRS-UMR 8089 et Universit\'{e} de Cergy-Pontoise\\
2 Avenue Adolphe Chauvin, F-95302, Cergy-Pontoise, France\\
E-mail: Thierry.Huillet@u-cergy.fr}
\maketitle

\begin{abstract}
We study a class of coalescents derived from a sampling procedure out of $N$
i.i.d. Pareto$\left( \alpha \right) $ random variables, normalized by their
sum, including $\beta -$size-biasing on total length effects ($\beta <\alpha 
$). Depending on the range of $\alpha ,$ we derive the large $N$ limit
coalescents structure, leading either to a discrete-time Poisson-Dirichlet$%
\left( \alpha ,-\beta \right) $ $\Xi -$coalescent ($\alpha \in \left[
0,1\right) $), or to a family of continuous-time Beta$\left( 2-\alpha
,\alpha -\beta \right) $ $\Lambda -$coalescents ($\alpha \in \left[
1,2\right) $), or to the Kingman coalescent ($\alpha \geq 2$). We indicate
that this class of coalescent processes (and their scaling limits) may be
viewed as the genealogical processes of some forward in time evolving
branching population models including selection effects. In such
constant-size population models, the reproduction step, which is based on a
fitness-dependent Poisson Point Process with scaling power-law$\left( \alpha
\right) $ intensity, is coupled to a selection step consisting of sorting
out the $N$ fittest individuals issued from the reproduction step.\newline

\textbf{Running title:} Pareto genealogies in a Poisson evolution model with
selection.\newline

\end{abstract}

\section{Introduction}

We first investigate discrete-time finite coalescents derived from sampling
from $N$ i.i.d. Pareto$\left( \alpha \right) $ random variables, normalized
by their sum, with $\alpha >0$. We include size-biasing on total length
effects involving a parameter $\beta <\alpha $. Depending on the range of $%
\alpha ,$ we derive the large $N$ limit coalescents structure: The case $%
\alpha \in \left[ 0,1\right) $ leads to a discrete-time Poisson-Dirichlet$%
\left( \alpha ,-\beta \right) $ $\Xi -$coalescent (with no time-scaling).
The case $\alpha =1$ gives rise to a continuous-time beta$\left( 1,1-\beta
\right) $ $\Lambda -$coalescent, involving a logarithmic time scaling $\log
N $. The case $\alpha \in \left( 1,2\right) $ leads to a continuous-time Beta%
$\left( 2-\alpha ,\alpha -\beta \right) $ $\Lambda -$coalescent, involving a
power-law time scaling according to $N^{\alpha -1}$. The range $\alpha \geq
2 $ gives rise to the standard Kingman coalescent (with time scaling $N$ if $%
\alpha >2$ and $N/\log N$ in the critical case $\alpha =2$). We give for
each case the exact speeds of convergence. We establish a loose link with a
Generalized Central Limit Theorem for stable random variables and we briefly
recall the main statistical features akin to general $\Lambda -$coalescents.

We indicate that the above special classes of coalescent processes (and
their scaling limits) may be viewed as the genealogical processes of some
forward in time evolving branching population models including selection
effects. These models are closely related in spirit to the additive
exponential model discussed in Brunet et al \cite{BDM}, \cite{BD}.

In the models we first consider, the size $N$ of the population is kept
constant over the generations. Each alive individual is assigned some
positive fitness $x>0$. In each generation and for each of the $N$ offspring
alive independently, the reproduction step is based on a fitness-dependent
Poisson Point Process (PPP) with scaling power-law$\left( \alpha \right) $
intensity; this procedure assigns new fitnesses to the (potentially
infinitely many) individuals of the next generation, in a multiplicative
way. We call $f\left( x\right) =x^{\alpha }$ the output fitness of an
individual with fitness $x$. The selection step then consists of sorting out
the $N$ fittest individuals issued from the reproduction step\footnote{%
This particular way of introducing selection in a randomly evolving
branching population with constant poulation size seems to appear first in 
\cite{BDMM}. It has nothing to do with the way selection is classically
introduced in population genetics; see \cite{Mar}, \cite{Ew} and \cite{RB}.}%
. The process is iterated independently over the subsequent generations. To
the first large $N$ approximation, the logarithm of the mean output fitness
within each generation $k$, scaled by the generation number $k$, is shown to
shift to the right at speed $v_{N}=\log \log N$ as $k\rightarrow \infty .$

While adopting a sampling point of view based on the intensity of the PPP to
compute the coalescence probabilities that some offspring is the one of a
parental individual with given fitness, it is shown that the genealogy of
the branching model with selection is in the domain of attraction of the beta%
$\left( 1,1-\beta \right) $ coalescent in the large $N$ limit
(Bolthausen-Sznitman coalescent if $\beta =0$). It is also shown that the
full class of the Pareto-coalescents discussed earlier can be obtained while
considering the PPP which is the output image of the original one, given by
the output map $f\left( x\right) =x^{\alpha }$ in fitness space$.$ In this
setup, the large $N$ limit computations of the coalescence and merging
probabilities are based not on the fitnesses but on the output deformed
fitnesses$.$

\section{Coalescents derived from Pareto-Sampling}

\subsection{Pareto sampling and coalescents}

Let $X_{1},..,X_{N}$ be $N$ i.i.d. Pareto$\left( \alpha \right) $ random
variables (rvs) with $\mathbf{P}\left( X_{1}>x\right) =x^{-\alpha }$, $%
\alpha >0$ and $x\geq 1.$ Let $F_{X_{1}}\left( x\right) =1-\mathbf{P}\left(
X_{1}>x\right) $ denote its probability distribution function (pdf)$.$ The
density of $X_{1}$ is $f_{X_{1}}\left( x\right) =\alpha x^{-\left( \alpha
+1\right) }.$ Let $S_{n}:=X_{n}/\sum_{1}^{N}X_{n}$, $n=1,...,N$ define a
random partition of the unit interval, upon normalizing the $X$s by their
sum. The $S_{n}$s are identically distributed but not independent of course
as they sum up to $1$; by doing so, the unit interval $\left[ 0,1\right] $
is thus broken into $N$ random pieces (subintervals or segments) of sizes $%
S_{n},$ $n=1,...,N.$

By sampling the $S$s, we mean that we draw independently $i$ uniform random
variables on the unit interval with $i\leq N$, looking at the subintervals
which are being hit in the process$.$ From this procedure, for instance, the
probability that the $i-$sample hits any one of the $S_{n}$s only once is 
\begin{equation}
P_{i,1}^{\left( N\right) }=\mathbf{E}\left( \sum_{n=1}^{N}S_{n}^{i}\right) =N%
\mathbf{E}\left( S_{1}^{i}\right) .  \label{P1}
\end{equation}
Let $\Sigma _{N}:=\sum_{n=1}^{N}X_{n}$ denote the partial sum of the $X$s.
For the values of $\beta $ for which $\mathbf{E}\left( \Sigma _{N}^{\beta
}\right) $ exists, we can size-bias the latter probability by the total
length $\Sigma _{N}$ and consider instead\footnote{%
We abusively use the same notation $P_{i,1}^{\left( N\right) }$ in the
size-biased setup as in (\ref{P1}) (corresponding to $\beta =0$), to avoid
overburden notations.} 
\begin{equation*}
P_{i,1}^{\left( N\right) }=\frac{\mathbf{E}\left( \Sigma _{N}^{\beta
}\sum_{n=1}^{N}S_{n}^{i}\right) }{\mathbf{E}\left( \Sigma _{N}^{\beta
}\right) }=\frac{N\mathbf{E}\left( \Sigma _{N}^{\beta }S_{1}^{i}\right) }{%
\mathbf{E}\left( \Sigma _{N}^{\beta }\right) }.
\end{equation*}
The latter event under consideration corresponds to an $i$ to $1$ merger of
some Markov coalescent process where $i$ particles are identified to a
single one (share the same ancestor) whenever the $i-$sample hits the same
subinterval of the unit partition. In this setup, $P_{i,1}^{\left( N\right) }
$ is therefore the entry $\left( i,1\right) $ of its one-step transition
matrix. The quantity $c_{N}:=P_{2,1}^{\left( N\right) },$ which is the
probability that two individuals chosen at random out of $N$ share the same
common ancestor, is called the coalescence probability.

Similarly we can define a $i$ to $j$ merger ($j\leq i$) by considering the
event that the $i$ particles hit any size$-j$ subset of the segments $S$
constituting the partition of unity. We get 
\begin{equation}
P_{i,j}^{\left( N\right) }=\binom{N}{j}\sum_{i_{1}+...+i_{j}=i}^{*}\binom{i}{%
i_{1}...i_{j}}\mathbf{E}\left( \prod_{l=1}^{j}S_{l}^{i_{l}}\right) =
\label{P2}
\end{equation}
\begin{equation*}
\binom{N}{j}\sum_{l=1}^{j}\left( -1\right) ^{j-l}\binom{j}{l}\mathbf{E}%
\left( \left( S_{1}+...+S_{l}\right) ^{i}\right) ,
\end{equation*}
where the star-sum in (\ref{P2}) runs over the $i_{l}\geq 1.$ The quantity $%
\mathbf{E}\left( \prod_{l=1}^{j}S_{l}^{i_{l}}\right) $ is the probability of
a $\left( i_{1},...,i_{j}\right) -$merger from $i$ to $j$. Using the same
abuse of notation, we shall also write the size-biased version of the latter
probability as 
\begin{equation}
P_{i,j}^{\left( N\right) }=\binom{N}{j}\sum_{i_{1}+...+i_{j}=i}^{*}\binom{i}{%
i_{1}...i_{j}}\frac{\mathbf{E}\left( \Sigma _{N}^{\beta
}\prod_{l=1}^{j}S_{l}^{i_{l}}\right) }{\mathbf{E}\left( \Sigma _{N}^{\beta
}\right) }.  \label{P2B}
\end{equation}
Unless stated otherwise, whenever we speak in the sequel of $P_{i,j}^{\left(
N\right) },$ we mean (\ref{P2B}) and not (\ref{P2}).

Clearly, $\sum_{j=1}^{i}P_{i,j}^{\left( N\right) }=1$ and so the matrix $%
P^{\left( N\right) }$ with entries $P_{i,j}^{\left( N\right) }$, $i=1,...,N$
and $j=1,...,i,$ is a $N\times N$ lower-triangular stochastic matrix
corresponding to some Markov discrete-time-$k$ coalescent (pure death)
process, say $x_{k}^{\left( N\right) }$, with finite state-space\footnote{%
A `true' coalescent process takes values in the set of equivalence relations
or partitions on $\left\{ 1,...,N\right\} $ and we rather deal here and
throughout with its block-counting counterpart.} and state $\left\{
1\right\} $ absorbing. Let us first investigate the expression of the
size-biased probability $P_{i,1}^{\left( N\right) },$ showing that its large 
$N$ estimate depends on the understanding of the $\beta -$moments of $\Sigma
_{N}.$

\begin{proposition}
When $-\infty <\beta <\alpha <2,$ the size-biased probability of an $i$ to $1
$ merger reads 
\begin{equation}
P_{i,1}^{\left( N\right) }=N\alpha \frac{\mathbf{E}\left( \Sigma
_{N-1}^{\beta -\alpha }\right) }{\mathbf{E}\left( \Sigma _{N}^{\beta
}\right) }\frac{\Gamma \left( i-\alpha \right) \Gamma \left( \alpha -\beta
\right) }{\Gamma \left( i-\beta \right) }\text{, }i\geq 2.  \label{P2s}
\end{equation}
\end{proposition}

\textbf{Proof:} We have $1/S_{1}=1+\Sigma _{N-1}^{\prime }/X_{1}$ where $%
\Sigma _{N-1}^{\prime }:=\sum_{n=2}^{N}X_{n}\overset{d}{=}%
\sum_{n=1}^{N-1}X_{n}=:\Sigma _{N-1}$. By conditioning on $X_{1},$ with $%
f_{\Sigma _{N-1}}$ the density of $\Sigma _{N-1}$, we get 
\begin{equation*}
P_{i,1}^{\left( N\right) }=\frac{N}{\mathbf{E}\left( \Sigma _{N}^{\beta
}\right) }\int_{1}^{\infty }dx\cdot f_{X_{1}}\left( x\right) x^{\beta
+1}\int_{1+\frac{N-1}{x}}^{\infty }z^{\beta -i}f_{\Sigma _{N-1}}\left(
x\left( z-1\right) \right) dz.
\end{equation*}
Reversing the integration and after two changes of variables 
\begin{eqnarray*}
P_{i,1}^{\left( N\right) } &=&\frac{N}{\mathbf{E}\left( \Sigma _{N}^{\beta
}\right) }\int_{1}^{\infty }dz\cdot z^{\beta -i}\int_{\frac{N-1}{z-1}%
}^{\infty }x^{\beta +1}f_{X_{1}}\left( x\right) f_{\Sigma _{N-1}}\left(
x\left( z-1\right) \right) dx \\
&=&\frac{N\alpha }{\mathbf{E}\left( \Sigma _{N}^{\beta }\right) }%
\int_{1}^{\infty }dz\cdot z^{\beta -i}\left( z-1\right) ^{\alpha -\beta
-1}\int_{N-1}^{\infty }s^{\beta -\alpha }f_{\Sigma _{N-1}}\left( s\right) ds
\\
&=&N\alpha \frac{\mathbf{E}\left( \Sigma _{N-1}^{\beta -\alpha }\right) }{%
\mathbf{E}\left( \Sigma _{N}^{\beta }\right) }\int_{0}^{1}u^{i-2}u^{1-\alpha
}\left( 1-u\right) ^{\alpha -\beta -1}du.
\end{eqnarray*}
When $\beta <\alpha <2,$ the latter integral is (upon adequate normalization
by a beta function term $B\left( 2-\alpha ,\alpha -\beta \right) $)
identified with the order $i-2$ moment of a beta($2-\alpha ,\alpha -\beta $)
rv. In this parameter range, for $i\geq 2$%
\begin{equation*}
\int_{0}^{1}u^{i-2}u^{1-\alpha }\left( 1-u\right) ^{\alpha -\beta -1}du=%
\frac{\Gamma \left( i-\alpha \right) \Gamma \left( \alpha -\beta \right) }{%
\Gamma \left( i-\beta \right) }=:B\left( i-\alpha ,\alpha -\beta \right)
\end{equation*}
are well-defined. $\diamond $\newline

The Pareto rv $X_{1}$ has power-law tails with index $\alpha .$ Therefore $%
\mathbf{E}\left( X_{1}^{\beta }\right) $ only exists when $\beta <\alpha ,$
with $\mathbf{E}\left( X_{1}^{\beta }\right) =\alpha /\left( \alpha -\beta
\right) .$ Clearly also, the tails of the sum $\Sigma _{N}$ obey $\mathbf{P}%
\left( \Sigma _{N}>s\right) \underset{s\rightarrow \infty }{\sim }N\mathbf{P}%
\left( X_{1}>s\right) =Ns^{-\alpha }.$ Because $N\mathbf{P}\left(
X_{1}>s\right) \underset{s\rightarrow \infty }{\sim }\mathbf{P}\left(
M_{N}>s\right) =1-\left( 1-s^{-\alpha }\right) ^{N}$ where $M_{N}=\max
\left( X_{1},..,X_{N}\right) $, this means that for large $s$, the event $%
\Sigma _{N}>s$ is essentially determined by the event $M_{N}>s$. Note that,
as a result, $\Sigma _{N}$ has the same tail index as $X_{1}$, indicating
that the $\beta -$moment of $\Sigma _{N}$ only exists for $\beta <\alpha .$%
\newline

We will show below from this, that large $N$ estimates of $\mathbf{E}\left(
\Sigma _{N}^{\beta }\right) $ can be obtained.

As a result, for example, based on (\ref{P2s}), whenever $1<\alpha <2$ and $%
\beta <\alpha ,$ it will be checked that 
\begin{equation*}
c_{N}:=P_{2,1}^{\left( N\right) }=N\alpha \frac{\mathbf{E}\left( \Sigma
_{N-1}^{\beta -\alpha }\right) }{\mathbf{E}\left( \Sigma _{N}^{\beta
}\right) }B\left( 2-\alpha ,\alpha -\beta \right) \propto N^{-\left( \alpha
-1\right) }\underset{N\rightarrow \infty }{\rightarrow }0.
\end{equation*}
We will show that, in that case, for each $i\geq 2,$ the limits $%
\lim_{N\rightarrow \infty }c_{N}^{-1}P_{i,1}^{\left( N\right) }$ exist, with 
\begin{equation*}
c_{N}^{-1}P_{i,1}^{\left( N\right) }\underset{N\rightarrow \infty }{%
\rightarrow }\frac{1}{B\left( 2-\alpha ,\alpha -\beta \right) }\frac{\Gamma
\left( \alpha -\beta \right) \Gamma \left( i-\alpha \right) }{\Gamma \left(
i-\beta \right) }=\int_{0}^{1}u^{i-2}\Lambda \left( du\right) >0,
\end{equation*}
where $\Lambda \left( du\right) $ has the density on $\left[ 0,1\right] $: $%
u^{1-\alpha }\left( 1-u\right) ^{\alpha -\beta -1}/B\left( 2-\alpha ,\alpha
-\beta \right) $, which is a beta($2-\alpha ,\alpha -\beta $) density.

Because in the large $N$ limit, simultaneous multiple collisions will be
seen to be negligible, we conclude, using similar arguments to the ones in 
\cite{Sa} and \cite{MS}, that, in the range $1<\alpha <2$ and $\beta <\alpha 
$, a time-scaled version of the finite discrete-time-$k$ coalescent $%
x_{k}^{\left( N\right) }$ arising from size-biased sampling out of Pareto$%
\left( \alpha \right) $ partitions converges weakly (as $N\rightarrow \infty 
$) to a continuous-time-$t$ $\Lambda -$coalescent $x_{t}$ where $\Lambda $
is a beta($2-\alpha ,\alpha -\beta $) probability measure on $\left[
0,1\right] ;$ namely $x_{t}\overset{d}{=}\lim_{N\rightarrow \infty
}x_{\left[ t/c_{N}\right] }^{\left( N\right) }$\footnote{%
Here, because two processes are involved, the symbol $\overset{d}{=}$ means
convergence of all the finite-dimensional distributions of $x_{\left[
t/c_{N}\right] }^{\left( N\right) },$ $t\geq 0$ to the ones of $x_{t},$ $%
t\geq 0$.}. In other words, if $t$ is continuous-time, the appropriate
scaling is $k\rightarrow \left[ t/c_{N}\right] $ (the integral part of $%
t/c_{N}$), showing that time $t$ should be measured in units of $%
N_{e}:=c_{N}^{-1}$ (the effective population size). The limiting process $%
x_{t}$ is a continuous-time pure death Markov process on $\Bbb{N}=\left\{
1,2,...\right\} $, absorbed at state $\left\{ 1\right\} $ and with
transition rates from $i$ to $j$ given by 
\begin{equation*}
c_{N}^{-1}P_{i,j}^{\left( N\right) }\underset{N\rightarrow \infty }{%
\rightarrow }\lambda _{i,j}:=\binom{i}{j-1}\int_{0}^{1}u^{i-j-1}\left(
1-u\right) ^{j-1}\Lambda \left( du\right) \text{, }1\leq j<i.
\end{equation*}
The rate terms $\lambda _{i,j}$ may also be written as 
\begin{equation*}
\lambda _{i,j}:=\sum_{l=0}^{j-1}\frac{\left( -1\right) ^{l}}{l!}\binom{i}{%
j-l-1}q_{l+i-j+1},
\end{equation*}
where $q_{l}:=\int_{0}^{1}u^{l-2}\Lambda \left( du\right) $ are the $l-$%
moments of $u^{-2}\Lambda \left( du\right) $ and also the rates of $l$ to $1$
mergers $\lambda _{l,1}.$

\subsection{Generalities on $\Lambda -$coalescents}

Whenever one gets a family of rates $\lambda _{i,j}$ as above for some
finite probability measure $\Lambda $ on $\left[ 0,1\right] ,$ one speaks of
continuous-time-$t$ $\Lambda -$coalescents (see \cite{Pit} for a precise
definition). These $\Lambda -$coalescents are non-increasing pure death
Markov processes, say $x_{t}$, on the state-space $\Bbb{N},$ and the $%
\lambda _{i,j}$s are the rates at which $i$ to $j<i$ mergers for $x_{t}$
occur. The state $\left\{ 1\right\} $ is absorbing. In such processes,
multiple collisions of any order (when $1\leq j<i$) can occur, but never
simultaneously. The total death rate at which some merger occurs, starting
from state $i$, is $\lambda _{i}:=\sum_{j=1}^{i-1}\lambda _{i,j}.$ One can
check that when $\Lambda =\delta _{0}$ (corresponding to the Kingman
coalescent), $\lambda _{i,j}\neq 0$ only when $j=i-1$ with $\lambda
_{i}=\lambda _{i,i-1}=\binom{i}{2}$. The general expression of $\lambda _{i}$
is 
\begin{equation*}
\lambda _{i}=\int_{\left[ 0,1\right] }u^{-2}\left( 1-\left( 1-u\right)
^{i}-iu\left( 1-u\right) ^{i-1}\right) \Lambda \left( du\right) .
\end{equation*}
When $\Lambda \left( \left\{ 0\right\} \right) =0$ (excluding thereby the
Kingman coalescent)$,$ the precise dynamics of $x_{t}$ when started at $i$
is given by $x_{0}=i$ and \cite{BB} 
\begin{equation}
x_{t}-x_{0}=-\int_{\left( 0,t\right] \times \left( 0,1\right] }\left(
B\left( x_{s_{-}},u\right) -1_{B\left( x_{s_{-}},u\right) >0}\right) 
\mathcal{N}\left( ds\times du\right)  \label{P2a}
\end{equation}
\begin{equation*}
=-\int_{\left( 0,t\right] \times \left( 0,1\right] }\left( B\left(
x_{s_{-}},u\right) -1\right) _{+}\mathcal{N}\left( ds\times du\right) .
\end{equation*}
Here, $x_{+}=\max \left( x,0\right) $, $\mathcal{N}$ is a random Poisson
measure on $\left[ 0,\infty \right) \times \left( 0,1\right] $ with
intensity $ds\times \frac{1}{u^{2}}\Lambda \left( du\right) $ and $B\left(
x_{s_{-}},u\right) \overset{d}{\sim }$ bin$\left( x_{s_{-}},u\right) $ is a
binomial rv with parameters $\left( x_{s_{-}},u\right) .$ As a result, with 
\begin{equation}
r\left( i\right) :=\int_{\left( 0,1\right] }\left( ui-1+\left( 1-u\right)
^{i}\right) u^{-2}\Lambda \left( du\right) ,\text{ }i>0,  \label{P2b}
\end{equation}
upon taking the expectation in (\ref{P2a}), it holds that 
\begin{equation*}
\mathbf{E}\left( dx_{t}\mid x_{t_{-}}\right) =-r\left( x_{t_{-}}\right) dt.
\end{equation*}
From this, the quantity $r\left( i\right) $ is the rate at which size $i$
blocks are being lost as time passes by. Clearly $r$ is also 
\begin{equation*}
r\left( i\right) =i\lambda _{i}-\sum_{j=1}^{i-1}j\lambda
_{i,j}=\sum_{j=1}^{i-1}\left( i-j\right) \lambda _{i,j}.
\end{equation*}
Consequently, the reciprocal function $1/r\left( i\right) $ of the rate $%
r\left( i\right) $ interprets as the expected time spent by $x_{t}$ in a
state with $i$ lineages and therefore $\sum_{j=2}^{x_{0}=i}1/r\left(
j\right) $ will give a large $i$ estimate of the expected time to the most
recent common ancestor (the height of the coalescent tree): 
\begin{equation*}
\tau _{i,1}:=\inf \left( t\in \Bbb{R}_{+}:x_{t}=1\mid x_{0}=i\right) .
\end{equation*}
There are lots of detailed studies in the literature (see a precise partial
list below) on other functionals of $x_{t}$ such as the total branch length $%
L_{i}$ of the $\Lambda -$coalescent, its total external branch length $%
L_{i}^{e}$, the length $l_{i}$ of its external branch (the time till first
collision of a branch chosen at random out of $i$), the number of collisions 
$C_{i}$ till time to most recent common ancestor,...\newline

All these functionals obey some distributional identities which prove useful
to obtain some insight on their limit laws as $i\rightarrow \infty $. Given $%
x_{0}=i$, all involve the number $U_{i}$ of singletons taking part in the
first collision occurring at time $T_{i}\overset{d}{\sim }$ exp$\left(
\lambda _{i}\right) ,$ giving $x_{T_{i}}-1=i-U_{i}$ singletons not
participating to the first collision (with $T_{i},$ $x_{T_{i}}$ independent)$%
.$ The quantity $U_{i}$ is important in itself because, due to $\mathbf{P}%
\left( U_{i}=j\right) =\lambda _{i,i-j+1}/\lambda _{i}$, $j=2,..,i$%
\begin{equation*}
r\left( i\right) =\lambda _{i}\left( \mathbf{E}\left( U_{i}\right) -1\right)
\end{equation*}
where 
\begin{equation*}
\mathbf{E}\left( U_{i}\right) =\frac{i}{\lambda _{i}}\int_{\left( 0,1\right]
}\frac{1-u-\left( 1-u\right) ^{i}}{u\left( 1-u\right) }\Lambda \left(
du\right) .
\end{equation*}
Provided $\Lambda $ has no atom at point $\left\{ 1\right\} ,$ the condition 
$\sum_{i=2}^{\infty }1/r\left( i\right) <\infty $ is the necessary and
sufficient condition for $x_{t}$ to come down from infinity, \cite{Sch3}.

Clearly indeed\footnote{$\overset{d}{=}$ means equality in distribution
between random variables.}, $\tau _{i,1}\overset{d}{=}T_{i}+\tau
_{x_{T_{i}},1}$, $L_{i}\overset{d}{=}iT_{i}+L_{x_{T_{i}}},$ $L_{i}^{e}%
\overset{d}{=}iT_{i}+L_{x_{T_{i}}-1}^{e},$ $C_{i}\overset{d}{=}%
1+C_{x_{T_{i}}}$ and $l_{i}\overset{d}{=}T_{i}+B_{i}l_{x_{T_{i}}-1}$ where $%
B_{i}$ is a Bernoulli rv, given by: $\mathbf{P}\left( B_{i}=1\mid
x_{T_{i}}\right) =\left( x_{T_{i}}-1\right) /i$ with $B_{i}l_{x_{T_{i}}-1}$
independent of $T_{i}$.\newline

Famous examples that we shall deal with in the sequel, include $\Lambda -$%
coalescents for which:

- $\Lambda \left( du\right) =B\left( a,b\right) u^{a-1}\left( 1-u\right)
^{b-1}1_{\left[ 0,1\right] }\left( u\right) du$, with $a,b>0$ and with $%
B\left( a,b\right) $ the beta function: we get the beta$\left( a,b\right) $\
coalescents.

- $\Lambda \left( du\right) =B\left( 2-\alpha ,\alpha \right) u^{1-\alpha
}\left( 1-u\right) ^{\alpha -1}1_{\left[ 0,1\right] }\left( u\right) du,$%
\thinspace ($\alpha \in \left( 0,1\right) \cup \left( 1,2\right) $); this is
the beta$\left( 2-\alpha ,\alpha \right) $ coalescent.

- (Lebesgue) $\Lambda \left( du\right) =1_{\left[ 0,1\right] }\left(
u\right) du:$\ this is the Bolthausen-Sznitman coalescent or beta$\left(
1,1\right) $\ coalescent.

- $\Lambda \left( du\right) =\delta _{0}:$ we get the Kingman coalescent
where only binary mergers can occur ($j=i-1$) one at a time, \cite{King}.

- $\Lambda \left( du\right) =\delta _{1}:$ we get the star-shaped coalescent
involving a single big collision. $\diamond $\newline

Using the above distributional identities, it can be shown for instance
that, for the Kingman coalescent and for large $i$, to leading order of
magnitude, rough estimates are: $\tau _{i,1}\sim 2\left( 1-1/i\right) $ \cite
{Tav}, $L_{i}\sim 2\log i$ \cite{DIMR}, $L_{i}^{e}\sim 2$ \cite{JK}, $%
C_{i}=i-1$ and $l_{i}\sim 1/i$ \cite{CNKR}$.$\newline

For the Bolthausen-Sznitman coalescent: $\tau _{i,1}\sim \log \log i$ \cite
{GM}, $L_{i}\sim i/\log i$ \cite{DIMR}, $C_{i}\sim i/\log i$ \cite{IM} and $%
l_{i}\sim 1/\log i$ \cite{FM}.\newline

For the beta$\left( 2-\alpha ,\alpha \right) $ coalescent with $0<\alpha <1$%
, $L_{i}\sim i$ \cite{M} and $L_{i}/L_{i}^{e}$ converges in probability to $%
1 $ \cite{Moh} and $l_{i}\sim O\left( 1\right) $ \cite{GIM}$.$ \newline

For the beta$\left( 2-\alpha ,\alpha \right) $ coalescent with $1<\alpha <2$%
, $L_{i}\sim 1/i^{\alpha -2}$ \cite{K}, $L_{i}^{e}\sim 1/i^{\alpha -2}$ \cite
{DY} and $l_{i}\sim 1/i^{\alpha -1}$ \cite{DFSY}$.$\newline

This information is important to grasp the general shape of the coalescent
trees in each case.

\section{GCLT for Pareto sums}

In this Section, we first sketch a loose connection of the large $N$
estimate of $c_{N}$ with the Generalized Central Limit Theorem for random
variables in the domain of attraction of stable laws, \cite{UZ}.

Let $\Sigma _{N}=\sum_{n=1}^{N}X_{n}$ be the partial sums of the i.i.d. $X$s
with Pareto$\left( \alpha \right) $ distributions, $\alpha >0,$ on $\left(
1,\infty \right) $. Let $S_{\alpha }$ be skewed $\alpha -$stable rvs with 
\begin{equation*}
\mathbf{E}\left( e^{-\lambda S_{\alpha }}\right) =e^{-\lambda ^{\alpha }}%
\text{ if }\alpha \in \left( 0,1\right) \text{, }\lambda \geq 0
\end{equation*}
the Laplace-Stieltjes transform (LST) of a one-sided $\alpha -$stable rv, 
\begin{equation*}
\mathbf{E}\left( e^{i\lambda S_{1}}\right) =e^{-\left| \lambda \right|
\left( 1+i\text{sign}\left( \lambda \right) \frac{2}{\pi }\log \left|
\lambda \right| \right) }
\end{equation*}
the characteristic function (c.f.) of a skewed $1-$stable Cauchy rv on $\Bbb{%
R}$ 
\begin{equation*}
\mathbf{E}\left( e^{i\lambda S_{\alpha }}\right) =e^{-\left| \lambda \right|
^{\alpha }\left( 1-i\text{sgn}\left( \lambda \right) \tan \left( \frac{\pi
\alpha }{2}\right) \right) }\text{ if }\alpha \in \left( 1,2\right) \text{, }%
\lambda \in \Bbb{R}
\end{equation*}
the characteristic function (c.f.) of a skewed $\alpha -$stable Cauchy rv on 
$\Bbb{R}$ and 
\begin{equation*}
\mathbf{E}\left( e^{i\lambda S_{\alpha }}\right) =e^{-\lambda ^{2}/2}\text{
if }\alpha \geq 2\text{, }\lambda \in \Bbb{R}
\end{equation*}
the c.f. of a standard normal rv on $\Bbb{R}.$

The following Generalized Central Limit Theorem (GCLT) then holds

\begin{theorem}
(\cite{UZ}, \cite{ZKS}): Let $\Sigma _{N}$ denote the partial sum sequence
of $N$ i.i.d. Pareto$\left( \alpha \right) $ random variables. Then, 
\begin{equation}
\frac{\Sigma _{N}-a_{N}}{b_{N}}\overset{d}{\underset{N\rightarrow \infty }{%
\rightarrow }}S_{\alpha }  \label{P3}
\end{equation}
where, with

\begin{eqnarray*}
C_{\alpha } &=&\left( \Gamma \left( 1-\alpha \right) \cos \left( \frac{\pi
\alpha }{2}\right) \right) ^{1/\alpha }\text{ if }\alpha \in \left(
0,2\right) \backslash \left\{ 1\right\} , \\
C_{1} &=&\frac{\pi }{2}, \\
C_{\alpha } &=&\left( \frac{\alpha }{\alpha -2}-\left( \frac{\alpha }{\alpha
-1}\right) ^{2}\right) ^{1/2}\text{ if }\alpha >2,
\end{eqnarray*}
$b_{N}$ is given by 
\begin{eqnarray*}
\left( i\right) \text{ }b_{N} &=&C_{\alpha }N^{\max \left( 1/\alpha
,1/2\right) }\text{ if }\alpha \neq 2 \\
\left( ii\right) \text{ }b_{N} &=&\left( N\log N\right) ^{1/2}\text{ if }%
\alpha =2
\end{eqnarray*}
and, with $\gamma $ the Euler constant and $\mathbf{E}\left( X_{1}\right)
=\mu :=\alpha /\left( \alpha -1\right) $, $a_{N}$ is given by 
\begin{eqnarray*}
\left( i\right) \text{ }a_{N} &=&0\text{ if }\alpha \in \left( 0,1\right)  \\
\left( ii\right) \text{ }a_{N} &=&\frac{N\pi ^{2}}{2}\int_{1}^{\infty }\sin
\left( \frac{2x}{\pi N}\right) dF_{X_{1}}\left( x\right) \sim N\log
N+N\left( 1-\gamma -\log \frac{2}{\pi }\right) \text{ if }\alpha =1 \\
\left( iii\right) \text{ }a_{N} &=&N\mu \text{ if }\alpha \in \left(
1,\infty \right) .
\end{eqnarray*}
\end{theorem}

When $\alpha \leq 1$, the characteristic values of $\Sigma _{N}$ can be
guessed to be what they are claimed to be ($N^{1/\alpha }$ if $\alpha <1$
and $N\log N$ if $\alpha =1$) while estimating $N\int_{1}^{m_{N}}xf_{X_{1}}%
\left( x\right) dx$ where $m_{N}$ is the mode of $M_{N}$ which is seen to
grow like $N^{1/\alpha }$ (see \cite{BG}) and similarly for the fluctuation
scaling term $b_{N}$ for $\alpha \leq 2$ (of order $N^{1/\alpha }$ if $%
\alpha <2$ and $\left( N\log N\right) ^{1/2}$ if $\alpha =2$).

From these rough estimates, we would conclude that with $c_{N}\propto N\frac{%
\mathbf{E}\left( \Sigma _{N-1}^{\beta -\alpha }\right) }{\mathbf{E}\left(
\Sigma _{N}^{\beta }\right) }$, up to the leading order in $N$%
\begin{equation*}
c_{N}\propto N\frac{\left( C_{\alpha }N^{1/\alpha }\right) ^{\beta -\alpha }%
}{\left( C_{\alpha }N^{1/\alpha }\right) ^{\beta }}=O\left( 1\right) \text{
if }\alpha \in \left( 0,1\right)
\end{equation*}
\begin{equation*}
c_{N}\propto N\frac{\left( N\log N\right) ^{\beta -1}}{\left( N\log N\right)
^{\beta }}\sim \frac{1}{\log N}\text{ if }\alpha =1
\end{equation*}
\begin{equation*}
c_{N}\propto N\frac{\left( \mu N\right) ^{\beta -\alpha }}{\left( \mu
N\right) ^{\beta }}\sim \mu ^{-\alpha }N^{-\left( \alpha -1\right) }\text{
if }\alpha \in \left( 1,2\right) .
\end{equation*}
Depending on the values of $\alpha $, we therefore anticipate\newline

$\bullet $ $\alpha \in \left( 0,1\right) $: Because in that case $c_{N}$ is
asymptotic to a constant, this suggests a limiting discrete-time coalescent$%
. $ We will show below that it is not a discrete-time $\Lambda -$coalescent,
rather it is a $\Xi -$coalescent of the Poisson-Dirichlet type with two
parameters $\left( \alpha ,-\beta \right) $. \newline
$\Xi -$coalescents were first introduced in \cite{MS} and further studied in 
\cite{Sch}. In sharp contrast with $\Lambda -$coalescents, multiple
collisions can occur simultaneously at the same transition time. In their
block-counting version, they are characterized by the set of numbers $\phi
_{j}\left( i_{1},..,i_{j}\right) $ defining the probabilities of an $\left(
i_{1},..,i_{j}\right) -$merger ($i_{1}+..+i_{j}=i$), resulting when the $\Xi
-$coalescent is discrete, in an $i$ to $j\leq i$ transition with probability 
$P_{i,j}=\frac{1}{j!}\sum_{i_{1}+...+i_{j}=i}^{*}\binom{i}{i_{1}...i_{j}}%
\phi _{j}\left( i_{1},..,i_{j}\right) $. The $\phi _{j}\left(
i_{1},..,i_{j}\right) $ can be written as 
\begin{equation*}
\phi _{j}\left( i_{1},..,i_{j}\right) =\int_{\Delta
_{j}}\prod_{l=1}^{j}u_{l}^{i_{l}-2}\Lambda _{j}\left(
du_{1},..,du_{j}\right) ,
\end{equation*}
for some finite measures $\Lambda _{j}$ with density on the $\left(
j+1\right) -$ simplex 
\begin{equation*}
\Delta _{j}=\left\{ \left( u_{1},..,u_{j}\right) \in \left[ 0,1\right]
^{j}:u_{1}+...+u_{j}\leq 1\right\} .
\end{equation*}
The set of measures $\Lambda _{j},$ $j\geq 1,$ (characterized by their
moments $\phi _{j}$), with values over the simplices $\Delta _{j},$
completely characterize the $\Xi -$coalescent, \cite{MS}.

In the simplest cases, with $\left\langle u,u\right\rangle :=\sum_{l\geq
1}u_{l}^{2}$, it is also convenient (see \cite{Sch}) to rewrite the $\phi
_{j}$s as 
\begin{equation*}
\phi _{j}\left( i_{1},..,i_{j}\right) =\int_{\Delta }\sum_{\underset{\text{%
all distinct}}{k_{1},...,k_{j}}}^{*}\prod_{l=1}^{j}u_{k_{l}}^{i_{l}}\frac{%
\Xi \left( du\right) }{\left\langle u,u\right\rangle },
\end{equation*}
where $\Delta =\left\{ u\mathbf{:=}\left( u_{1},..,u_{l},..\right)
:u_{1}\geq ..\geq u_{l}\geq ..\geq 0:\sum_{l\geq 1}u_{l}\leq 1\right\} $ and 
$\Xi $ a finite measure concentrated on the subset $\Delta ^{*}$ of $\Delta $
consisting of those $u$s exactly summing to $1.$ Letting $\nu \left(
du\right) :=\Xi \left( du\right) /\left\langle u,u\right\rangle ,$ the
measure $\nu $ on the infinite simplex $\Delta $ is such that $\nu \left(
\Delta \right) <\infty .$\newline

$\bullet $ $\alpha =1$: A $\Lambda -$coalescent with logarithmic effective
population size where $\Lambda $ is a beta($1,1-\beta $) probability measure
with $\beta <1$ as in \cite{BD} (reducing to the Bolthausen-Sznitman
coalescent when avoiding size-biasing corresponding to $\beta =0$).

$\bullet $ $\alpha \in \left( 1,2\right) $: A $\Lambda -$coalescent with $%
N_{e}\propto N^{\alpha -1}$ where $\Lambda $ is a beta($2-\alpha ,\alpha
-\beta $) measure with $\beta <\alpha .$ In the latter case, $\beta =0$
leads to the standard beta($2-\alpha ,\alpha $) coalescent.\newline

$\bullet $ $\alpha \geq 2:$ In the range $\alpha >2$, (respectively $\alpha
=2$), $N_{e}\propto N$ (respectively $N_{e}\propto N/\log N$) and the
obtained scaled limiting coalescent is a Kingman coalescent, as a result of $%
\Lambda $ approaching $\delta _{0}$. In this latter case, only binary
mergers occur in the large $N$ limit.\newline

\textbf{Remark:} The Kingman coalescent also occurs when dealing with
sampling from some alternative random partition. For instance, would
sampling be defined from a random partition of unity given by $%
S_{n}:=X_{n}/\sum_{1}^{N}X_{n}$, $n=1,...,N$ where $X_{1}$ now obeys the
following gamma$\left( \theta \right) $ density: $f_{X_{1}}\left( x\right)
=\Gamma \left( \theta \right) ^{-1}x^{\theta -1}e^{-x}$, $\theta ,x>0,$ then
the law of $\Sigma _{N}=X_{1}+...+X_{N}$ is $f_{\Sigma _{N}}\left( x\right)
=\Gamma \left( N\theta \right) ^{-1}x^{N\theta -1}e^{-x}$, independent of $%
S_{1}$ and 
\begin{eqnarray*}
P_{i,1}^{\left( N\right) } &=&\frac{\mathbf{E}\left( \Sigma _{N}^{\beta
}S_{n}^{i}\right) }{\mathbf{E}\left( \Sigma _{N}^{\beta }\right) }=\frac{N%
\mathbf{E}\left( \Sigma _{N}^{\beta }S_{1}^{i}\right) }{\mathbf{E}\left(
\Sigma _{N}^{\beta }\right) }=N\mathbf{E}\left( S_{1}^{i}\right) \\
&=&N\frac{\Gamma \left( N\theta \right) }{\Gamma \left( \theta \right) }%
\frac{\Gamma \left( i+\theta \right) }{\Gamma \left( N\theta +i\right) }\sim
N^{-\left( i-1\right) }\frac{\Gamma \left( i+\theta \right) }{\Gamma \left(
\theta \right) }\theta ^{-i}.
\end{eqnarray*}
Thus $c_{N}=P_{2,1}^{\left( N\right) }=\frac{1}{N}\frac{1+\theta }{\theta }%
\underset{N\rightarrow \infty }{\rightarrow }0$ together with $%
d_{N}=P_{3,1}^{\left( N\right) }=\frac{1}{N^{2}}\frac{\left( 1+\theta
\right) \left( 2+\theta \right) }{\theta ^{2}}$. Because triple mergers are
asymptotically negligible compared to binary ones ($d_{N}/c_{N}\rightarrow 0$%
), the time-scaled limiting coalescent using $c_{N}=\frac{1}{N}\frac{%
1+\theta }{\theta }\propto N^{-1}$ is a Kingman coalescent. The prefactor $%
\frac{1+\theta }{\theta }$ appearing in front of $c_{N}$ is the ratio $\rho
/\mu ^{2}$ where $\rho :=\mathbf{E}\left( X_{1}^{2}\right) =\theta \left(
\theta +1\right) $ and $\mu :=\mathbf{E}\left( X_{1}\right) =\theta $ and $%
c_{N}$ is independent of $\beta .$

\section{Large $N$ asymptotic estimation of $\mathbf{E}\left( \Sigma
_{N}^{\beta }\right) $ and consequences}

In this Section, we compute the asymptotic behavior of the $\beta -$moments
of $\Sigma _{N}$ in the cases $\alpha \in \left( 0,1\right) ,$ $\alpha \in
\left( 1,2\right) $ and $\alpha =1$, and $\alpha \geq 2,$ making the
previous conclusions based on the GCLT consistent\footnote{%
The technique we use is inspired from the one used in \cite{BD} in a
particular case. The author is indebted to B. Derrida for pointing this out
to him.}$.$\newline

$\bullet $ We start with the case $\alpha \in \left( 0,1\right) .$

\begin{theorem}
When $\alpha \in \left( 0,1\right) ,$ as $N\rightarrow \infty ,$ $x_{k}%
\overset{d}{=}\lim_{N\rightarrow \infty }x_{k}^{\left( N\right) }$ exists
and is a discrete-time Poisson-Dirichlet$\left( \alpha ,-\beta \right) $ $%
\Xi -$coalescent$.$
\end{theorem}

\textbf{Proof:} For the values of $\beta $ for which it makes sense, we have 
\begin{equation*}
\mathbf{E}\left( \Sigma _{N}^{\beta }\right) =\frac{1}{\Gamma \left( -\beta
\right) }\int_{0}^{\infty }d\lambda \cdot \lambda ^{-\beta -1}\mathbf{E}%
\left( e^{-\lambda \Sigma _{N}^{{}}}\right) =\frac{1}{\Gamma \left( -\beta
\right) }\int_{0}^{\infty }d\lambda \cdot \lambda ^{-\beta -1}\mathbf{E}%
\left( e^{-\lambda X_{1}}\right) ^{N}.
\end{equation*}
When $N$ is large, only the small $\lambda $ approximation of $\mathbf{E}%
\left( e^{-\lambda X_{1}}\right) $ to the latter integral contributes. For
small $\lambda $, we have 
\begin{equation}
\mathbf{E}\left( e^{-\lambda X_{1}}\right) =\alpha \int_{1}^{\infty
}x^{-\left( \alpha +1\right) }e^{-\lambda x}dx=1-\alpha \int_{1}^{\infty
}x^{-\left( \alpha +1\right) }\left( 1-e^{-\lambda x}\right) dx\sim
\label{PI}
\end{equation}
\begin{equation*}
1-\alpha \int_{0}^{\infty }x^{-\left( \alpha +1\right) }\left( 1-e^{-\lambda
x}\right) dx\sim 1-\Gamma \left( 1-\alpha \right) \lambda ^{\alpha }\sim
e^{-\Gamma \left( 1-\alpha \right) \lambda ^{\alpha }}.
\end{equation*}
Note that, when $\lambda $ is small 
\begin{equation*}
\mathbf{E}\left( X_{1}^{i}e^{-\lambda X_{1}}\right) =\left( -1\right) ^{i}%
\frac{d^{i}}{d\lambda ^{i}}\mathbf{E}\left( e^{-\lambda X_{1}}\right) \sim
\alpha \Gamma \left( i-\alpha \right) \lambda ^{\alpha -i}.
\end{equation*}
Thus, 
\begin{equation*}
\mathbf{E}\left( \Sigma _{N}^{\beta }\right) \sim \frac{1}{\Gamma \left(
-\beta \right) }\int_{0}^{\infty }d\lambda \cdot \lambda ^{-\beta
-1}e^{-N\Gamma \left( 1-\alpha \right) \lambda ^{\alpha }}.
\end{equation*}
With the change of variables $u=N\Gamma \left( 1-\alpha \right) \lambda
^{\alpha }$, with $\beta <\alpha $, we get 
\begin{equation*}
\mathbf{E}\left( \Sigma _{N}^{\beta }\right) \sim \frac{N^{\beta /\alpha }}{%
\alpha \Gamma \left( -\beta \right) }\Gamma \left( 1-\alpha \right) ^{\beta
/\alpha }\int_{0}^{\infty }du\cdot u^{-\beta /\alpha -1}e^{-u}=N^{\beta
/\alpha }\frac{\Gamma \left( 1-\alpha \right) ^{\beta /\alpha }\Gamma \left(
1-\beta /\alpha \right) }{\Gamma \left( 1-\beta \right) }.
\end{equation*}
Finally, using the above large $N$ estimate of $\mathbf{E}\left( \Sigma
_{N}^{\beta }\right) $, the identity $\Gamma \left( x+1\right) =x\Gamma
\left( x\right) $ and (\ref{P2s}) with $i=2,$ when $\beta <\alpha <1,$ we
obtain 
\begin{eqnarray*}
c_{N} &=&N\alpha \frac{\mathbf{E}\left( \Sigma _{N-1}^{\beta -\alpha
}\right) }{\mathbf{E}\left( \Sigma _{N}^{\beta }\right) }\frac{\Gamma \left(
2-\alpha \right) \Gamma \left( \alpha -\beta \right) }{\Gamma \left( 2-\beta
\right) }\sim \alpha \frac{\Gamma \left( 1-\beta \right) \left( 1-\frac{%
\beta }{\alpha }\right) }{\Gamma \left( 1-\left( \beta -\alpha \right)
\right) \Gamma \left( 1-\alpha \right) }\frac{\Gamma \left( 2-\alpha \right)
\Gamma \left( \alpha -\beta \right) }{\Gamma \left( 2-\beta \right) } \\
&=&\frac{\Gamma \left( 1-\beta \right) }{\Gamma \left( \alpha -\beta \right)
\Gamma \left( 1-\alpha \right) }\frac{\Gamma \left( 2-\alpha \right) \Gamma
\left( \alpha -\beta \right) }{\Gamma \left( 2-\beta \right) }=\frac{%
1-\alpha }{1-\beta }=:c_{\infty }.
\end{eqnarray*}
Thus the coalescence probability $c_{N}$ converges to $c_{\infty }\in \left(
0,1\right) .$ When $\alpha \in \left( 0,1\right) $, using again (\ref{P2s}),
we obtain more generally 
\begin{equation*}
P_{i,1}^{\left( N\right) }\underset{N\rightarrow \infty }{\rightarrow }%
P_{i,1}=\frac{\Gamma \left( 1-\beta \right) }{\Gamma \left( 1-\alpha \right) 
}\frac{\Gamma \left( i-\alpha \right) }{\Gamma \left( i-\beta \right) },
\end{equation*}
which are the probabilities to merge all $i$ particles in one step in the
limiting discrete-time coalescent.

To derive the full transition probabilities of the limiting discrete-time-$k$
coalescent $x_{k}\overset{d}{=}\lim_{N\rightarrow \infty }x_{k}^{\left(
N\right) }$, recalling that 
\begin{equation*}
P_{i,j}^{\left( N\right) }=\binom{N}{j}\sum_{i_{1}+...+i_{j}=i}^{*}\binom{i}{%
i_{1}...i_{j}}\frac{\mathbf{E}\left( \Sigma _{N}^{\beta
}\prod_{l=1}^{j}S_{l}^{i_{l}}\right) }{\mathbf{E}\left( \Sigma _{N}^{\beta
}\right) },
\end{equation*}
we use 
\begin{equation}
\mathbf{E}\left( \Sigma _{N}^{\beta }\prod_{l=1}^{j}S_{l}^{i_{l}}\right) =%
\frac{1}{\Gamma \left( i-\beta \right) }\int_{0}^{\infty }d\lambda \cdot
\lambda ^{i-\beta -1}\prod_{l=1}^{j}\mathbf{E}\left(
X_{l}^{i_{l}}e^{-\lambda X_{l}}\right) \mathbf{E}\left( e^{-\lambda
X_{1}}\right) ^{N-j}  \label{P0}
\end{equation}
\begin{equation*}
\sim \frac{\alpha ^{j}}{\Gamma \left( i-\beta \right) }\prod_{l=1}^{j}\Gamma
\left( i_{l}-\alpha \right) \int_{0}^{\infty }d\lambda \cdot \lambda
^{i-\beta -1}\lambda ^{\alpha j-i}e^{-N\Gamma \left( 1-\alpha \right)
\lambda ^{\alpha }}.
\end{equation*}
Performing again the change of variables $u=N\Gamma \left( 1-\alpha \right)
\lambda ^{\alpha }$, with $\beta <\alpha $, we get 
\begin{equation*}
\mathbf{E}\left( \Sigma _{N}^{\beta }\prod_{l=1}^{j}S_{l}^{i_{l}}\right)
\sim \alpha ^{j-1}N^{\beta /\alpha -j}\Gamma \left( 1-\alpha \right) ^{\beta
/\alpha -1}\frac{\Gamma \left( j-\beta /\alpha \right) }{\Gamma \left(
i-\beta \right) }\prod_{l=1}^{j}\Gamma \left( i_{l}-\alpha \right) .
\end{equation*}
Using $\binom{N}{j}\sim N^{j}/j!$ for large $N$, we finally get $%
P_{i,j}^{\left( N\right) }\rightarrow P_{i,j}$ where ($1\leq j\leq i$) 
\begin{equation}
P_{i,j}=\frac{i!}{j!}\alpha ^{j-1}\frac{\Gamma \left( 1-\beta \right) }{%
\Gamma \left( 1-\beta /\alpha \right) }\frac{\Gamma \left( j-\beta /\alpha
\right) }{\Gamma \left( i-\beta \right) }\sum_{i_{1}+...+i_{j}=i}^{*}%
\prod_{l=1}^{j}\frac{\Gamma \left( i_{l}-\alpha \right) }{\Gamma \left(
1-\alpha \right) i_{l}!}.  \label{P4}
\end{equation}
These are the full transition probabilities of the limiting discrete-time-$k$
coalescent $x_{k}$ in the regime $\alpha \in \left( 0,1\right) $ and $\beta
<\alpha $ (satisfying $\sum_{j=1}^{i}P_{i,j}=1$)$.$ Note that $P_{i,1}$ are
the probabilities obtained previously and that the diagonal terms (the
eigenvalues of $P$) read 
\begin{equation*}
P_{i,i}=\alpha ^{i-1}\frac{\Gamma \left( 1-\beta \right) }{\Gamma \left(
1-\beta /\alpha \right) }\frac{\Gamma \left( i-\beta /\alpha \right) }{%
\Gamma \left( i-\beta \right) }.
\end{equation*}
Clearly this discrete-time coalescent is not a discrete $\Lambda -$%
coalescent as simultaneous multiple collisions can occur (it is a $\Xi -$%
coalescent). Clearly $P_{i,j}$ is also 
\begin{equation*}
P_{i,j}=\frac{1}{j!}\sum_{i_{1}+...+i_{j}=i}^{*}\binom{i}{i_{1}...i_{j}}\phi
_{j}\left( i_{1},..,i_{j}\right) ,
\end{equation*}
where the $\phi _{j}\left( i_{1},..,i_{j}\right) $s define the probabilities
of a $\left( i_{1},..,i_{j}\right) -$merger ($i_{1}+..+i_{j}=i$). These $%
\phi _{j}\left( i_{1},..,i_{j}\right) ,$ which can be read from (\ref{P4}),
may be written under the alternative form 
\begin{eqnarray*}
\phi _{j}\left( i_{1},..,i_{j}\right) &:&=\alpha ^{j-1}\frac{\Gamma \left(
1-\beta \right) }{\Gamma \left( 1-\beta /\alpha \right) }\frac{\Gamma \left(
j-\beta /\alpha \right) }{\Gamma \left( i-\beta \right) }\prod_{l=1}^{j}%
\frac{\Gamma \left( i_{l}-\alpha \right) }{\Gamma \left( 1-\alpha \right) }
\\
&=&c_{j,\alpha ,\beta }\frac{\Gamma \left( \alpha j-\beta \right) }{\Gamma
\left( i-\beta \right) }\prod_{l=1}^{j}\Gamma \left( i_{l}-\alpha \right) ,
\end{eqnarray*}
where 
\begin{equation*}
c_{j,\alpha ,\beta }:=\prod_{l=1}^{j}\frac{\Gamma \left( \left( l-1\right)
\alpha +1-\beta \right) }{\Gamma \left( 1-\alpha \right) \Gamma \left(
l\alpha -\beta \right) }.
\end{equation*}
Thus (\ref{P4}) is also 
\begin{equation}
P_{i,j}=c_{j,\alpha ,\beta }\frac{i!}{j!}\frac{\Gamma \left( \alpha j-\beta
\right) }{\Gamma \left( i-\beta \right) }\sum_{i_{1}+...+i_{j}=i}^{*}%
\prod_{l=1}^{j}\frac{\Gamma \left( i_{l}-\alpha \right) }{i_{l}!}.
\label{P4a}
\end{equation}
Defining the finite Dirichlet measures $\Lambda _{j}$ with density on the $%
\left( j+1\right) -$ simplex $\Delta _{j}$ given by:

\begin{equation*}
\Lambda _{j}\left( du_{1},..,du_{j}\right) =\frac{c_{j,\alpha ,\beta }}{%
c_{\infty }}\prod_{l=1}^{j}\left( u_{l}^{1-\alpha }du_{l}\right) \left(
1-\sum_{l=1}^{j}u_{l}\right) ^{\alpha j-\beta -1},
\end{equation*}
we get 
\begin{equation*}
\phi _{j}\left( i_{1},..,i_{j}\right) =c_{\infty }\int_{\Delta
_{j}}\prod_{l=1}^{j}u_{l}^{i_{l}-2}\Lambda _{j}\left(
du_{1},..,du_{j}\right) .
\end{equation*}
The set of finite Dirichlet measures $\Lambda _{j}$ with parameters 
\begin{equation*}
\left( \theta _{1}=2-\alpha ,...,\theta _{j}=2-\alpha ,\theta _{j+1}=\alpha
j-\beta \right)
\end{equation*}
on the simplices $\Delta _{j}$ completely characterize this limiting
discrete-time coalescent. Note that $\Lambda _{1}$ is a beta$\left( 2-\alpha
,\alpha -\beta \right) $ probability measure.

One may rewrite the $\phi _{j}$s as (see \cite{Sch} and \cite{Moh}) 
\begin{equation*}
\phi _{j}\left( i_{1},..,i_{j}\right) =\int_{\Delta }\sum_{\underset{\text{%
all distinct}}{k_{1},...,k_{j}}}^{*}\prod_{l=1}^{j}u_{k_{l}}^{i_{l}}\frac{%
\Xi \left( du\right) }{\left\langle u,u\right\rangle },
\end{equation*}
where $\Delta =\left\{ u\mathbf{:=}\left( u_{1},..,u_{l},..\right)
:u_{1}\geq ..\geq u_{l}\geq ..\geq 0:\sum_{l\geq 1}u_{l}\leq 1\right\} $ and 
$\Xi $ a measure on $\Delta .$ Letting $\nu \left( du\right) :=\Xi \left(
du\right) /\left\langle u,u\right\rangle ,$ the measure $\nu $ on the
infinite simplex $\Delta $ can be identified (see \cite{Moh}) to the
two-parameter Poisson-Dirichlet$\left( \alpha ,-\beta \right) $ measure,
with $\alpha \in \left[ 0,1\right) $ and $\beta <\alpha .$ It holds that $%
\nu \left( \Delta \right) =1.$ Poisson-Dirichlet measures enjoy many
remarkable properties including a stick-breaking property, \cite{PitYor}. $%
\diamond $\newline

\textbf{Remarks:}

$\left( i\right) $ From (\ref{P4a}), the limiting situation $\alpha =0$ also
makes sense, leading to the one-parameter Poisson-Dirichlet$\left( 0,-\beta
\right) $ measure with $\beta <\alpha =0$. In this case, from (\ref{P4a}) 
\begin{equation*}
P_{i,j}=\frac{\left( -\beta \right) ^{j}\Gamma \left( -\beta \right) }{%
\Gamma \left( i-\beta \right) }s_{i,j},
\end{equation*}
where $s_{i,j}:=\frac{i!}{j!}$ $\sum_{i_{1}+...+i_{j}=i}^{*}\prod_{l=1}^{j}%
\frac{1}{i_{l}}$ are the absolute first kind Stirling numbers.\newline

$\left( ii\right) $ Finally, avoiding size-biasing ($\beta =0$) gives rise
to the discrete Poisson-Dirichlet$\left( \alpha ,0\right) $ coalescent with
one-parameter $\alpha \in \left( 0,1\right) ,$ appearing in \cite{Sch2}.%
\newline

$\bullet $ The case $\alpha \in \left( 1,2\right) .$

\begin{lemma}
When $\alpha \in \left( 1,2\right) ,$ with $\mu :=\frac{\alpha }{\alpha -1},$
for all $\beta <\alpha $%
\begin{equation}
c_{N}:=P_{2,1}^{\left( N\right) }\sim \alpha \mu ^{-\alpha }B\left( 2-\alpha
,\alpha -\beta \right) N^{-\left( \alpha -1\right) }\underset{N\rightarrow
\infty }{\rightarrow }0.  \label{P4b}
\end{equation}
\end{lemma}

\textbf{Proof:} In this case, defining $a:=\alpha -1\in \left( 0,1\right) ,$
after an integration by parts and using the previous estimate (\ref{PI})
substituting $a$ to $\alpha $%
\begin{equation*}
\mathbf{E}\left( e^{-\lambda X_{1}}\right) =e^{-\lambda }-\frac{\lambda }{a}%
a\int_{1}^{\infty }x^{-\left( a+1\right) }e^{-\lambda x}dx\sim e^{-\lambda }-%
\frac{\lambda }{a}\left( 1-\Gamma \left( 1-a\right) \lambda ^{a}\right) .
\end{equation*}
Thus, for small $\lambda ,$ with $\mu :=\frac{\alpha }{\alpha -1}$%
\begin{equation}
\mathbf{E}\left( e^{-\lambda X_{1}}\right) \sim 1-\lambda \mu -\Gamma \left(
1-\alpha \right) \lambda ^{\alpha }\sim e^{-\lambda \mu }.  \label{PII}
\end{equation}
Thus, consistently with the GCLT approach, to the dominant order 
\begin{equation*}
\mathbf{E}\left( \Sigma _{N}^{\beta }\right) =\frac{1}{\Gamma \left( -\beta
\right) }\int_{0}^{\infty }d\lambda \cdot \lambda ^{-\beta -1}\mathbf{E}%
\left( e^{-\lambda X_{1}}\right) ^{N}\sim \left( N\mu \right) ^{\beta }\text{%
,}
\end{equation*}
so that 
\begin{equation*}
c_{N}:=N\alpha \frac{\mathbf{E}\left( \Sigma _{N-1}^{\beta -\alpha }\right) 
}{\mathbf{E}\left( \Sigma _{N}^{\beta }\right) }B\left( 2-\alpha ,\alpha
-\beta \right) \sim \alpha \mu ^{-\alpha }B\left( 2-\alpha ,\alpha -\beta
\right) N^{-\left( \alpha -1\right) }.\text{ }\diamond
\end{equation*}

This suggests that

\begin{proposition}
When $\alpha \in \left( 1,2\right) ,$ upon scaling time using an effective
population size $N_{e}=c_{N}^{-1}$ (with $c_{N}$ as in (\ref{P4b}))$,$ we
obtain a limiting continuous-time-$t$ $\Lambda -$coalescent: $x_{t}\overset{d%
}{=}\lim_{N\rightarrow \infty }x_{\left[ t/c_{N}\right] }^{\left( N\right) },
$ with $\Lambda $ a beta($2-\alpha ,\alpha -\beta $) probability measure, $%
\beta <\alpha .$ If $\beta =0,$ we get the standard beta($2-\alpha ,\alpha $%
) coalescent.
\end{proposition}

\textbf{Proof:} To confirm this point, we will first evaluate a large $N$\
estimate of $P_{i,1}^{\left( N\right) }$\ defined in (\ref{P2}), in the
range $\alpha \in \left( 1,2\right) .$\ Using the above small $\lambda $\
estimate of $\mathbf{E}\left( e^{-\lambda X_{1}}\right) $ in the parameter
range under concern, for $i\geq 2,$\ we get 
\begin{equation*}
\mathbf{E}\left( X_{1}^{i}e^{-\lambda X_{1}}\right) =\left( -1\right) ^{i}%
\frac{d^{i}}{d\lambda ^{i}}\mathbf{E}\left( e^{-\lambda X_{1}}\right) \sim
\alpha \Gamma \left( i-\alpha \right) \lambda ^{\alpha -i}.
\end{equation*}
Thus, using (\ref{P0}) 
\begin{equation*}
\mathbf{E}\left( \Sigma _{N}^{\beta }S_{1}^{i}\right) =\frac{1}{\Gamma
\left( i-\beta \right) }\int_{0}^{\infty }d\lambda \cdot \lambda ^{i-\beta
-1}\mathbf{E}\left( X_{1}^{i}e^{-\lambda X_{1}}\right) \mathbf{E}\left(
e^{-\lambda X_{1}}\right) ^{N-1}
\end{equation*}
\begin{equation*}
\sim \alpha \frac{\Gamma \left( i-\alpha \right) }{\Gamma \left( i-\beta
\right) }\int_{0}^{\infty }d\lambda \cdot \lambda ^{i-\beta -1}\lambda
^{\alpha -i}e^{-\mu \lambda N}=\alpha \frac{\Gamma \left( i-\alpha \right) }{%
\Gamma \left( i-\beta \right) }\frac{\Gamma \left( \alpha -\beta \right) }{%
\left( \mu N\right) ^{\alpha -\beta }}.
\end{equation*}
Finally, we obtain 
\begin{equation*}
P_{i,1}^{\left( N\right) }\sim N^{-\left( \alpha -1\right) }\frac{\alpha }{%
\mu ^{\alpha }}\frac{\Gamma \left( \alpha -\beta \right) \Gamma \left(
i-\alpha \right) }{\Gamma \left( i-\beta \right) },
\end{equation*}
showing that, with $c_{N}=P_{2,1}^{\left( N\right) }\propto N^{-\left(
\alpha -1\right) }\rightarrow 0$, for each $i,$ $\lim_{N\rightarrow \infty
}c_{N}^{-1}P_{i,1}^{\left( N\right) }$ exist and are strictly positive
constants$.$ More precisely,

\begin{equation*}
c_{N}^{-1}P_{i,1}^{\left( N\right) }\underset{N\rightarrow \infty }{%
\rightarrow }\phi _{1}\left( i\right) =\frac{1}{B\left( 2-\alpha ,\alpha
-\beta \right) }\frac{\Gamma \left( \alpha -\beta \right) \Gamma \left(
i-\alpha \right) }{\Gamma \left( i-\beta \right) }=\int_{0}^{1}u^{i-2}%
\Lambda _{1}\left( du\right) ,
\end{equation*}
where $\Lambda _{1}=\Lambda $ is a beta($2-\alpha ,\alpha -\beta $)
probability measure.

To deal with the higher order terms, with $m:=\left\{ l\in \left\{
1,...,j\right\} :i_{l}\geq 2\right\} $, assuming $1\leq j<i,$ let us write (%
\ref{P2B}) as 
\begin{equation*}
P_{i,j}^{\left( N\right) }=\binom{N}{j}\sum_{m=1}^{j}\binom{j}{m}%
\sum_{i_{1}+...+i_{m}=i-j+m}^{**}\binom{i}{i_{1}...i_{m}}\frac{\mathbf{E}%
\left( \Sigma _{N}^{\beta
}\prod_{l=1}^{m}S_{l}^{i_{l}}\prod_{l=m+1}^{j}S_{l}\right) }{\mathbf{E}%
\left( \Sigma _{N}^{\beta }\right) },
\end{equation*}
where the double-star sum is now over the $i_{l},$ $l=1,...,m$ satisfying $%
i_{l}\geq 2.$ Proceeding similarly with higher order terms, it clearly holds
that for all $j\geq 2$ and $i_{l}\geq 1,$ $l=1,...,j,$ satisfying $i_{l}\geq
2$ for at least two $l$ in the list$,$ 
\begin{equation*}
\lim_{N\rightarrow \infty }c_{N}^{-1}\mathbf{E}\left( \Sigma _{N}^{\beta
}\prod_{l=1}^{j}S_{l}^{i_{l}}\right) /\mathbf{E}\left( \Sigma _{N}^{\beta
}\right) =0,
\end{equation*}
so that simultaneous multiple collisions cannot occur in the limit (actually
the contribution of these simultaneous multiple collisions terms in $%
P_{i,j}^{\left( N\right) }$ is of order $O\left( c_{N}^{m}\right) $). In the
latter expression of $P_{i,j}^{\left( N\right) }$ therefore, only the term
corresponding to $m=1$ will contribute to the $O\left( c_{N}\right) -$order.
Because, $\mathbf{E}\left( X_{1}e^{-\lambda X_{1}}\right) =-\frac{d}{%
d\lambda }\mathbf{E}\left( e^{-\lambda X_{1}}\right) \sim \mu $ and $\mathbf{%
E}\left( X_{1}^{i_{1}}e^{-\lambda X_{1}}\right) \sim \alpha \Gamma \left(
i_{1}-\alpha \right) \lambda ^{\alpha -i_{1}}$ with $i_{1}=i-j+1\geq 2$,
using (\ref{P0}), we get 
\begin{equation*}
\frac{\mathbf{E}\left( \Sigma _{N}^{\beta
}S_{1}^{i_{1}}\prod_{l=2}^{j}S_{l}\right) }{\mathbf{E}\left( \Sigma
_{N}^{\beta }\right) }=\int_{0}^{\infty }d\lambda \cdot \lambda ^{i-\beta -1}%
\mathbf{E}\left( X_{1}^{i_{1}}e^{-\lambda X_{1}}\right) \frac{\prod_{l=2}^{j}%
\mathbf{E}\left( X_{l}e^{-\lambda X_{l}}\right) \mathbf{E}\left( e^{-\lambda
X_{1}}\right) ^{N-j}}{\mathbf{E}\left( \Sigma _{N}^{\beta }\right) \Gamma
\left( i-\beta \right) }
\end{equation*}
\begin{equation*}
\sim \frac{\alpha \mu ^{j-1}\Gamma \left( i_{1}-\alpha \right) }{\mathbf{E}%
\left( \Sigma _{N}^{\beta }\right) \Gamma \left( i-\beta \right) }%
\int_{0}^{\infty }d\lambda \cdot \lambda ^{i-\beta -1}\lambda ^{\alpha
-i_{1}}e^{-N\mu \lambda }=\frac{\alpha \mu ^{j-1}\Gamma \left( i_{1}-\alpha
\right) }{\mathbf{E}\left( \Sigma _{N}^{\beta }\right) \Gamma \left( i-\beta
\right) }\frac{\Gamma \left( \alpha -\beta +j-1\right) }{\left( N\mu \right)
^{\alpha -\beta +j-1}}.
\end{equation*}
This shows, using $\binom{N}{j}\sim N^{j}/j!,$ (\ref{P4b}) and after some
elementary algebra, that for $j=1,...,i-1,$%
\begin{eqnarray*}
c_{N}^{-1}P_{i,j}^{\left( N\right) }\underset{N\rightarrow \infty }{%
\rightarrow }\lambda _{i,j} &=&\binom{i}{j-1}\frac{B\left( i-j+1-\alpha
,\alpha -\beta +j-1\right) }{B\left( 2-\alpha ,\alpha -\beta \right) } \\
&=&\binom{i}{j-1}\int_{0}^{1}u^{i-j-1}\left( 1-u\right) ^{j-1}\Lambda \left(
du\right) ,
\end{eqnarray*}
where $\Lambda =\Lambda _{1}\overset{d}{\sim }$ beta($2-\alpha ,\alpha
-\beta $). As a result, $x_{t}\overset{d}{=}\lim_{N\rightarrow \infty
}x_{\left[ t/c_{N}\right] }^{\left( N\right) }$ is a continuous-time-$t$
pure death coalescent process on $\Bbb{N}$ with infinitesimal transition
rates $\lambda _{i,j}.$ $\diamond $\newline

$\bullet $ The case $\alpha =1.$ We view it as a limiting case of the
previous analysis deriving from (\ref{PII}) when $\alpha \rightarrow 1^{+}.$

\begin{proposition}
When $\alpha =1,$ upon scaling time using a logarithmic effective population
size $N_{e}=c_{N}^{-1}\sim \log N,$ the limiting process $x_{t}\overset{d}{=}%
\lim_{N\rightarrow \infty }x_{\left[ t/c_{N}\right] }^{\left( N\right) }$
exists and is a continuous-time-$t$ $\Lambda -$coalescent with $\Lambda $ a
beta($1,1-\beta $) probability measure, $\beta <1.$ If $\beta =0,$ we get
the standard Bolthausen-Sznitman coalescent with $\Lambda $ uniform.
\end{proposition}

\textbf{Proof:} Put indeed $\alpha =1+\varepsilon $ in the previous small-$%
\lambda $ estimate (\ref{PII}) of $\mathbf{E}\left( e^{-\lambda
X_{1}}\right) $ in the parameter range $\alpha \in \left( 1,2\right) $, with 
$\varepsilon >0$ small. Then $\mu \sim 1+1/\varepsilon $ and because $\Gamma
\left( 1-\alpha \right) =\Gamma \left( -\varepsilon \right) \sim
_{0^{+}}-1/\varepsilon $, 
\begin{equation}
\mathbf{E}\left( e^{-\lambda X_{1}}\right) \sim 1-\lambda -\frac{1}{%
\varepsilon }\left( \lambda -\lambda ^{1+\varepsilon }\right) \sim 1-\lambda
+\lambda \log \lambda =:I\left( \lambda \right) .  \label{PIII}
\end{equation}
Thus, 
\begin{equation*}
\mathbf{E}\left( \Sigma _{N}^{\beta }\right) \sim \frac{1}{\Gamma \left(
-\beta \right) }\int_{0}^{\infty }d\lambda \cdot \lambda ^{-\beta -1}I\left(
\lambda \right) ^{N}.
\end{equation*}
The leading contribution of $I\left( \lambda \right) ^{N}$ is, when $\lambda 
$ is small, of order $1/\left( N\log N\right) $. Putting $\lambda =u/\left(
N\log N\right) $%
\begin{equation*}
I\left( \lambda \right) ^{N}\sim e^{-u}\left( 1+\frac{u}{\log N}\left( \log
u-\log \log N-1\right) \right) .
\end{equation*}
Thus, with $\beta <1$ 
\begin{equation*}
\mathbf{E}\left( \Sigma _{N}^{\beta }\right) \sim \frac{\left( N\log
N\right) ^{\beta }}{\Gamma \left( -\beta \right) }\int_{0}^{\infty }du\cdot
u^{-\beta -1}e^{-u}\left( 1+\frac{u}{\log N}\left( \log u-\log \log
N-1\right) \right)
\end{equation*}
\begin{equation}
\sim \left( N\log N\right) ^{\beta }\left( 1+\frac{\beta }{\log N}\left(
\psi \left( 1-\beta \right) -\log \log N-1\right) \right) ,  \label{P5}
\end{equation}
where $\psi \left( x\right) =\Gamma ^{\prime }\left( x\right) /\Gamma \left(
x\right) $ is the digamma function. In the latter estimate, we used that
differentiating $\int_{0}^{\infty }du\cdot u^{-\beta }u^{\theta
}e^{-u}=\Gamma \left( 1-\beta +\theta \right) $ with respect to the extra
parameter $\theta $ and then putting $\theta =0$ gives $\int_{0}^{\infty
}du\cdot u^{-\beta }\log \left( u\right) e^{-u}=\Gamma ^{\prime }\left(
1-\beta \right) $.

Finally, consistently with the GCLT approach 
\begin{equation*}
c_{N}=N\frac{\mathbf{E}\left( \Sigma _{N-1}^{\beta -1}\right) }{\mathbf{E}%
\left( \Sigma _{N}^{\beta }\right) }\sim \frac{1}{\log N}.
\end{equation*}
Clearly, with $\Lambda $ a beta($1,1-\beta $) probability measure and $\beta
<\alpha =1$%
\begin{equation*}
c_{N}^{-1}P_{i,1}^{\left( N\right) }\underset{N\rightarrow \infty }{%
\rightarrow }\int_{0}^{1}u^{i-2}\Lambda \left( du\right) =\frac{\Gamma
\left( 2-\beta \right) \Gamma \left( i-1\right) }{\Gamma \left( i-\beta
\right) }=\lambda _{i,1}.
\end{equation*}
Using similar arguments on higher order terms, showing that simultaneous
multiple collisions do not contribute in the limit, one can easily show 
\begin{equation}
c_{N}^{-1}P_{i,j}^{\left( N\right) }\underset{N\rightarrow \infty }{%
\rightarrow }\int_{0}^{1}u^{i-j-1}\left( 1-u\right) ^{j-1}\Lambda \left(
du\right)  \label{P5a}
\end{equation}
\begin{equation*}
=\frac{B\left( i-j,j-\beta \right) }{B\left( 1,1-\beta \right) }=:\lambda
_{i,j},\text{ }j=1,...,i-1,
\end{equation*}
where $\Lambda =\Lambda _{1}\overset{d}{\sim }$ beta($1,1-\beta $). This
confirms that, when $\alpha =1,$ upon scaling time using an effective
logarithmic population size $N_{e}=c_{N}^{-1}\sim \log N,$ we get a limiting
continuous-time-$t$ $\Lambda -$coalescent $x_{t}\overset{d}{=}%
\lim_{N\rightarrow \infty }x_{\left[ t/c_{N}\right] }^{\left( N\right) },$
where $\Lambda $ is a beta($1,1-\beta $) probability measure, $\beta <1.$ $%
\diamond $\newline

$\bullet $ The case $\alpha >2.$

\begin{proposition}
When $\alpha >2,$ for any value of $\beta $, with $\mu :=\frac{\alpha }{%
\alpha -1}>0$ and $\rho :=\frac{\alpha }{\alpha -2}>0$, upon scaling time
using a linear effective population size $N_{e}=c_{N}^{-1}=N\mu ^{2}/\rho ,$
the limiting process $x_{t}\overset{d}{=}\lim_{N\rightarrow \infty
}x_{\left[ t/c_{N}\right] }^{\left( N\right) }$ exists and is the
continuous-time-$t$ Kingman coalescent$.$
\end{proposition}

\textbf{Proof:} If $\alpha >2$, $\Sigma _{N}$\ is in the domain of
attraction of the normal law. As a result, for small $\lambda ,$ with $\mu :=%
\mathbf{E}\left( X_{1}\right) =\frac{\alpha }{\alpha -1}$ and $\rho :=%
\mathbf{E}\left( X_{1}^{2}\right) =\frac{\alpha }{\alpha -2}$%
\begin{equation}
\mathbf{E}\left( e^{-\lambda X_{1}}\right) \sim 1-\lambda \mu +\frac{1}{2}%
\rho \lambda ^{2}\sim e^{-\lambda \mu }.  \label{PIV}
\end{equation}
From (\ref{PIV}), for small $\lambda ,$ $\mathbf{E}\left( X_{1}e^{-\lambda
X_{1}}\right) =-\frac{d}{d\lambda }\mathbf{E}\left( e^{-\lambda
X_{1}}\right) \sim \mu $ and $\mathbf{E}\left( X_{1}^{2}e^{-\lambda
X_{1}}\right) =\frac{d^{2}}{d\lambda ^{2}}\mathbf{E}\left( e^{-\lambda
X_{1}}\right) \sim \rho $.

We have $c_{N}:=P_{2,1}^{\left( N\right) }=N\frac{\mathbf{E}\left( \Sigma
_{N}^{\beta }S_{1}^{2}\right) }{\mathbf{E}\left( \Sigma _{N}^{\beta }\right) 
}$ with, from (\ref{P0}) 
\begin{eqnarray*}
\mathbf{E}\left( \Sigma _{N}^{\beta }S_{1}^{2}\right) &\sim &\frac{\rho }{%
\Gamma \left( 2-\beta \right) }\int_{0}^{\infty }d\lambda \cdot \lambda
^{1-\beta }e^{-\lambda \mu N}=\rho \left( \mu N\right) ^{\beta -2} \\
\mathbf{E}\left( \Sigma _{N}^{\beta }\right) &\sim &\frac{1}{\Gamma \left(
-\beta \right) }\int_{0}^{\infty }d\lambda \cdot \lambda ^{-\beta
-1}e^{-\lambda \mu N}=\left( \mu N\right) ^{\beta }.
\end{eqnarray*}
Thus, 
\begin{equation*}
c_{N}=N\frac{\mathbf{E}\left( \Sigma _{N}^{\beta }S_{1}^{2}\right) }{\mathbf{%
E}\left( \Sigma _{N}^{\beta }\right) }\sim \frac{\rho }{\mu ^{2}}N^{-1}
\end{equation*}
goes to $0$ as claimed. As in Proposition $5$, simultaneous multiple
collisions cannot contribute in the limit. Only the term in the expression
of $P_{i,j}^{\left( N\right) }$ corresponding to $m=1$ will contribute to
the $O\left( c_{N}\right) -$order and so we need to focus on a collision
with a single $i_{1}\geq 2$.

Now, from (\ref{PIV}), $\mathbf{E}\left( X_{1}e^{-\lambda X_{1}}\right) =-%
\frac{d}{d\lambda }\mathbf{E}\left( e^{-\lambda X_{1}}\right) \sim \mu $ and 
$\mathbf{E}\left( X_{1}^{i_{1}}e^{-\lambda X_{1}}\right) \sim \rho \delta
_{i_{1},2}$ with $i_{1}=i-j+1\geq 2$. Proceeding again as in Proposition $5$
therefore, only effective transitions from $i$ to $j<i$ with $%
j=i-i_{1}+1=i-1 $ are seen in the limit, corresponding to the binary mergers
of a Kingman coalescent$.$ $\diamond $\newline

$\bullet $ $\alpha =2.$ To complete the picture, it remains to study the
limiting critical case $\alpha =2.$

\begin{lemma}
When $\alpha =2$, for all values of $\beta $, the coalescence probability
goes to $0$ as $N\rightarrow \infty $ like 
\begin{equation*}
c_{N}\sim \frac{1}{2}\frac{\log N}{N}.
\end{equation*}
\end{lemma}

\textbf{Proof:} We view the case $\alpha =2$ as a limiting case of (\ref{PIV}%
) as $\alpha \rightarrow 2^{+}.$\ When $\lambda $\ is small and for $\alpha
>2$, we have 
\begin{equation*}
\mathbf{E}\left( e^{-\lambda X_{1}}\right) \sim 1-\lambda \mu +\frac{1}{2}%
\rho \lambda ^{2}-\Gamma \left( 1-\alpha \right) \lambda ^{\alpha }.
\end{equation*}
Putting $\alpha =2+\varepsilon ,$\ for $\varepsilon >0$\ small, using $%
\Gamma \left( -1-\varepsilon \right) \sim _{0^{+}}1/\varepsilon $, with $%
\rho \sim 2/\varepsilon ,$\ when $\lambda $\ is small, we have 
\begin{equation}
\mathbf{E}\left( e^{-\lambda X_{1}}\right) \sim 1-2\lambda +\frac{1}{%
\varepsilon }\lambda ^{2}-\frac{1}{\varepsilon }\lambda ^{2+\varepsilon
}\sim 1-2\lambda -\lambda ^{2}\log \lambda \sim e^{-2\lambda }.  \label{PV}
\end{equation}
Because $\mathbf{E}\left( X_{1}e^{-\lambda X_{1}}\right) =-\frac{d}{d\lambda 
}\mathbf{E}\left( e^{-\lambda X_{1}}\right) \sim 2$\ and $\mathbf{E}\left(
X_{1}^{2}e^{-\lambda X_{1}}\right) =\frac{d^{2}}{d\lambda ^{2}}\mathbf{E}%
\left( e^{-\lambda X_{1}}\right) \sim -2\log \lambda -3$, we have 
\begin{eqnarray*}
\mathbf{E}\left( \Sigma _{N}^{\beta }S_{1}^{2}\right) &\sim &\frac{1}{\Gamma
\left( 2-\beta \right) }\int_{0}^{\infty }d\lambda \cdot \lambda ^{1-\beta }%
\mathbf{E}\left( X_{1}^{2}e^{-\lambda X_{1}}\right) \mathbf{E}\left(
e^{-\lambda X_{1}}\right) ^{N-1}\sim \\
&\sim &\frac{-1}{\Gamma \left( 2-\beta \right) }\int_{0}^{\infty }d\lambda
\cdot \lambda ^{1-\beta }\left( 2\log \lambda +3\right) e^{-2N\lambda } \\
\mathbf{E}\left( \Sigma _{N}^{\beta }\right) &\sim &\frac{1}{\Gamma \left(
-\beta \right) }\int_{0}^{\infty }d\lambda \cdot \lambda ^{-\beta
-1}e^{-2\lambda N}=\left( 2N\right) ^{\beta }.
\end{eqnarray*}
The Euler integral with a logarithmic term inside appearing in the
expression of $\mathbf{E}\left( \Sigma _{N}^{\beta }S_{1}^{2}\right) $ can
be obtained while taking the derivative of $\int_{0}^{\infty }d\lambda \cdot
\lambda ^{1-\beta }\lambda ^{\theta }e^{-2N\lambda }$ with respect to the
extra parameter $\theta $ and then putting $\theta =0$ in the obtained
expression. Observing therefore 
\begin{equation*}
\int_{0}^{\infty }d\lambda \cdot \lambda ^{1-\beta }\log \left( \lambda
\right) e^{-2N\lambda }=\Gamma ^{\prime }\left( 2-\beta \right) \left(
2N\right) ^{-\left( 2-\beta \right) }-\Gamma \left( 2-\beta \right) \left(
2N\right) ^{-\left( 2-\beta \right) }\log \left( 2N\right) ,
\end{equation*}
to the dominant order in $N$, 
\begin{equation*}
\mathbf{E}\left( \Sigma _{N}^{\beta }S_{1}^{2}\right) \sim 2\left( 2N\right)
^{-\left( 2-\beta \right) }\log N,
\end{equation*}
leading to 
\begin{equation*}
c_{N}:=P_{2,1}^{\left( N\right) }=N\frac{\mathbf{E}\left( \Sigma _{N}^{\beta
}S_{1}^{2}\right) }{\mathbf{E}\left( \Sigma _{N}^{\beta }\right) }\sim \frac{%
1}{2}\frac{\log N}{N}.\text{ }\diamond
\end{equation*}

\textbf{Remark:} A similar scaling behavior for\textbf{\ }$c_{N}$ was
recently obtained in Theorem $2.4$ of \cite{HM2}, dealing with coalescents
arising from compound Poisson discrete reproduction models, in the critical
case.

\begin{proposition}
When $\alpha =2$, upon scaling time using an effective population size $%
N_{e}=c_{N}^{-1}=\left( 2N\right) /\log N,$ the limiting process $x_{t}%
\overset{d}{=}\lim_{N\rightarrow \infty }x_{\left[ t/c_{N}\right] }^{\left(
N\right) }$ is the continuous-time-$t$ Kingman coalescent.
\end{proposition}

\textbf{Proof:} Using (\ref{PV}), for small $\lambda $, we have 
\begin{equation*}
\mathbf{E}\left( X_{1}^{i}e^{-\lambda X_{1}}\right) =\left( -1\right) ^{i}%
\frac{d^{i}}{d\lambda ^{i}}\mathbf{E}\left( e^{-\lambda X_{1}}\right) \sim
2\Gamma \left( i-2\right) \lambda ^{-\left( i-2\right) },\text{ }i\geq 3,
\end{equation*}
to which one should add $\mathbf{E}\left( X_{1}e^{-\lambda X_{1}}\right)
\sim 2$ and $\mathbf{E}\left( X_{1}^{2}e^{-\lambda X_{1}}\right) \sim
-\left( 2\log \lambda +3\right) .$ For $i\geq 3$, we get 
\begin{eqnarray*}
\mathbf{E}\left( \Sigma _{N}^{\beta }S_{1}^{i}\right) &=&\frac{1}{\Gamma
\left( i-\beta \right) }\int_{0}^{\infty }d\lambda \cdot \lambda ^{i-\beta
-1}\mathbf{E}\left( X_{1}^{i}e^{-\lambda X_{1}}\right) \mathbf{E}\left(
e^{-\lambda X_{1}}\right) ^{N-1} \\
&\sim &\frac{2\Gamma \left( i-2\right) }{\Gamma \left( i-\beta \right) }%
\int_{0}^{\infty }d\lambda \cdot \lambda ^{1-\beta }e^{-2N\lambda }=2\frac{%
B\left( i-2,2-\beta \right) }{\left( 2N\right) ^{2-\beta }}.
\end{eqnarray*}
Thus, for all $i\geq 3$, as $N\rightarrow \infty $%
\begin{equation*}
c_{N}^{-1}P_{i,1}^{\left( N\right) }=Nc_{N}^{-1}\frac{\mathbf{E}\left(
\Sigma _{N}^{\beta }S_{1}^{i}\right) }{\mathbf{E}\left( \Sigma _{N}^{\beta
}\right) }\sim \frac{1}{\log N}B\left( i-2,2-\beta \right) \rightarrow 0.
\end{equation*}
Therefore, due to the extra factor $\log N$ appearing in $c_{N}$, the
transitions from $i\geq 3$ to $1$ cannot be seen in the limit, nor (for the
same reason) the transitions involving a single multiple collision with $%
i_{1}\geq 3$, nor transitions involving simultaneous multiple collisions of
any order. In fact, only the events involving a single multiple collision
with $i_{1}=i-j+1=2$ (corresponding to transitions from $i$ to $j=i-1$) will
contribute in the limit. Indeed, 
\begin{equation*}
\frac{\mathbf{E}\left( \Sigma _{N}^{\beta
}S_{1}^{2}\prod_{l=2}^{j}S_{l}\right) }{\mathbf{E}\left( \Sigma _{N}^{\beta
}\right) }=\int_{0}^{\infty }d\lambda \cdot \lambda ^{i-\beta -1}\mathbf{E}%
\left( X_{1}^{2}e^{-\lambda X_{1}}\right) \frac{\prod_{l=2}^{j}\mathbf{E}%
\left( X_{l}e^{-\lambda X_{l}}\right) \mathbf{E}\left( e^{-\lambda
X_{1}}\right) ^{N-j}}{\mathbf{E}\left( \Sigma _{N}^{\beta }\right) \Gamma
\left( i-\beta \right) }
\end{equation*}
\begin{equation*}
\sim \frac{-2^{j-1}}{\mathbf{E}\left( \Sigma _{N}^{\beta }\right) \Gamma
\left( i-\beta \right) }\int_{0}^{\infty }d\lambda \cdot \lambda ^{i-\beta
-1}\left( 2\log \lambda +3\right) e^{-2N\lambda }\sim \frac{1}{2}N^{-i}\log
N.
\end{equation*}
This shows, using 
\begin{equation*}
P_{i,j}^{\left( N\right) }\sim \binom{N}{j}\binom{j}{1}\binom{i}{2\text{ }%
1...1}\frac{\mathbf{E}\left( \Sigma _{N}^{\beta
}S_{1}^{2}\prod_{l=2}^{j}S_{l}\right) }{\mathbf{E}\left( \Sigma _{N}^{\beta
}\right) }
\end{equation*}
with $j=i-1,$ $\binom{N}{j}\sim N^{j}/j!,$ and after some elementary
algebra, that

\begin{equation*}
c_{N}^{-1}P_{i,i-1}^{\left( N\right) }\underset{N\rightarrow \infty }{%
\rightarrow }\lambda _{i,i-1}=\binom{i}{2}.\text{ }\diamond
\end{equation*}

\textbf{Remarks:}

$\left( i\right) $ Results of a similar flavor can be found in
Schweinsberg's work \cite{Sch2}. However, his model and techniques are
different from ours because he considers large$-N$ limiting coalescents
obtained while sampling without replacement from a discrete super-critical
Galton-Watson branching process, assuming the reproduction law of each
offspring to exhibit power-law Zipf tails of index $\alpha $ ($a$ in his
notations). Note that there is no parameter $\beta $ in the construction 
\cite{Sch2}. In the same spirit, results can also be found in
Huillet-M\"{o}hle \cite{HM} where, following \cite{EW}, $\Lambda -$%
coalescents are obtained as scaling limits of discrete extended Moran
models, the skewed reproduction law of which displaying occasional extreme
events with one individual allowed to produce a large amount of offspring.
Whenever the reproduction law displays systematic extreme events, discrete
coalescents were even shown to emerge in the large $N-$limit, but in this
Moran context, they are only $\Lambda -$coalescents with multiple but no
simultaneous collisions, \cite{HM}. This contrasts with the occurrence in
our present work of discrete Poisson-Dirichlet $\Xi -$coalescents. In
contrast also with our current work where coalescents are derived from a
sampling procedure in the continuum, the coalescents considered in \cite
{Sch2} and \cite{HM} were built from discrete reproduction laws at fixed
population size $N$ and looking at their scaling limits $N\rightarrow \infty 
$.

We refer to these two works and to \cite{BD} for additional background on $%
\Lambda -$coalescent processes adapted to our purposes.\newline

$\left( ii\right) $ When $\alpha \geq 2$, neither the (large $N$ estimate of
the) scaling constant $c_{N}$ nor the limiting (Kingman) coalescent depend
on the bias parameter $\beta .$\newline

$\left( iii\right) $ We observe that the probability $c_{N}:=P_{2,1}^{\left(
N\right) }$ that two individuals chosen at random share the same common
ancestor, defining the time-scale to derive the large-$N$ limits of the
Pareto-coalescents all have the same large-$N$ order of magnitude as the
length $l_{N}$ of an external branch chosen at random in the limiting $%
\Lambda -$coalescents $x_{t}$ ($\alpha \in \left[ 1,2\right) $) or $\Xi -$%
coalescent $x_{k}$ ($\alpha \in \left( 0,1\right) $), or Kingman coalescent $%
x_{t}$ ($\alpha >2$), started at $x_{0}=N$. This curious fact is unexplained
so far.

\section{Forward in time selection model and genealogies}

In this Section, in the spirit of \cite{BD}, we indicate that the coalescent
processes just discussed may be viewed as the genealogical processes of some
forward in time evolving branching population models with selection. As it
is often the case in population genetics, the process we are interested in
is in the class of branching processes conditioned on having a fixed
population size over each generation, in the spirit of \cite{KMG}.

\subsection{A Poisson-point process model with selection}

Start with $N$ individuals at generation $t=0$ and assume that each
individual has an initial fitness $x_{n}\left( 0\right) >0$, $n=1,...,N.$

To describe the state of the population at the next generation, assume first
that, independently of one another, each individual potentially generates an
infinite number of offspring along a Poisson point process (PPP) with
intensity (or occupation) density 
\begin{equation}
\pi _{x_{n}\left( 0\right) }\left( x\right) =-\overline{\pi }_{x_{n}\left(
0\right) }^{\prime }\left( x\right) =\alpha x_{n}\left( 0\right) ^{\alpha
}x^{-\left( \alpha +1\right) },\footnote{%
The $^{\prime }$ symbol indicates derivative with respect to $x.$}
\label{B1}
\end{equation}
where $\overline{\pi }_{x_{n}\left( 0\right) }\left( x\right) :=\left(
x/x_{n}\left( 0\right) \right) ^{-\alpha }$, $x>0,$ $\alpha >0,$ $n=1,..,N.$
We observe that, with $\overline{\pi }\left( x\right) :=x^{-\alpha }$%
\begin{equation*}
\overline{\pi }_{x_{n}\left( 0\right) }\left( x\right) =\frac{\overline{\pi }%
\left( x\right) }{\overline{\pi }\left( x_{n}\left( 0\right) \right) }
\end{equation*}
and if we let $\pi \left( x\right) =-\overline{\pi }^{\prime }\left(
x\right) =\alpha x^{-\left( \alpha +1\right) },$ then $\pi _{x_{n}\left(
0\right) }\left( x\right) =x_{n}\left( 0\right) ^{\alpha }\pi \left(
x\right) .$

So the fitnesses of the offspring of each individual is generated according
to a PPP depending on the fitness of its parent. From these simple
assumptions, and as conventional wisdom suggests, we get:

\begin{proposition}
In a PPP model for fitness-dependent offspring reproduction with occupation
density (\ref{B1}), the fitnesses of the offspring of some parental
individual with fitness $x_{n}\left( 0\right) $ are $x_{n}\left( 0\right) $
times the fitnesses of the offspring of some canonical individual with unit
fitness: In this sense, the larger the fitness $x_{n}\left( 0\right) $ of
some individual is, the more he will, proportionally, produce offspring with
large fitness.
\end{proposition}

\textbf{Proof:} Let $\left( \tau _{n};n\geq 1\right) $ be the points of a
standard homogeneous Poisson point process (PPP) on the half-line with rate $%
1$, and let $\overline{\pi }^{-1}\left( s\right) =s^{-1/\alpha }$ be the
decreasing inverse of $\overline{\pi }.$ Then, with $\overline{\pi }%
_{x_{n}\left( 0\right) }^{-1}\left( s\right) =x_{n}\left( 0\right) \cdot
s^{-1/\alpha }$, $\left( \overline{\pi }_{x_{n}\left( 0\right) }^{-1}\left(
\tau _{n}\right) ;n\geq 1\right) $ are the (ordered) points of the offspring
PPP on the positive half-line (or here the fitness space) with occupation
density $\pi _{x_{n}\left( 0\right) }$. So, the fitter the individuals, the
fitter their offspring, proportionally to the parental fitness. $\diamond $%
\newline

Let us also briefly emphasize that, avoiding the fitness dependence on the
parent of the PPP, $\left( \overline{\pi }^{-1}\left( \tau _{n}\right)
;n\geq 1\right) $ are just the (ordered) points of a PPP on the half-line
with occupation density $\pi $. Whenever, as in our case study here, the
rate function $\pi $ is not integrable up to $x=0,$ there are infinitely
many such points, with $0$ as an accumulation point whereas there is of
course a finite Poisson (with mean $\overline{\pi }\left( \varepsilon
\right) $) number of them above some threshold $\varepsilon >0$.

It is well-known that when $\alpha \in \left( 0,1\right) $, with $\overline{%
\pi }\left( x\right) :=x^{-\alpha }$ and $\overline{\pi }^{-1}\left(
s\right) =s^{-1/\alpha }$, the positive cumulative rv 
\begin{equation}
\chi \overset{d}{=}\sum_{n\geq 1}\overline{\pi }^{-1}\left( \tau _{n}\right)
\label{Pcum}
\end{equation}
is a one-sided $\alpha -$stable rv on $\left( 0,\infty \right) $ with LST $%
\mathbf{E}\left( e^{-\lambda \chi }\right) =e^{-\kappa \lambda ^{\alpha }}$, 
$\kappa =\Gamma \left( 1-\alpha \right) >0$, $\lambda \geq 0.$

When $\alpha \in \left( 1,2\right) $, the law of $\chi $ is the one of a
positive Lamperti rv with LST $\mathbf{E}\left( e^{-\lambda \chi }\right)
=e^{-c\lambda +\kappa \lambda ^{\alpha }}$, $\kappa =-\Gamma \left( 1-\alpha
\right) >0,$ $c:=\mathbf{E}\left( \chi \right) >0$, \cite{Lam2}. When $%
\alpha =1$, the law of $\chi $ is the one of a positive Neveu rv with LST $%
\mathbf{E}\left( e^{-\lambda \chi }\right) =e^{\lambda \log \lambda }$ (see 
\cite{Neveu}); the latter may be viewed as a Lamperti rv in the limit $%
\alpha \rightarrow 1^{+},$ \cite{H2}$.$\newline

Typically indeed, by the L\'{e}vy-Khintchine formula, $\pi \left( x\right)
dx $ stands for the L\'{e}vy measure for the non-negative jumps of $\chi $
which is the value at time $t=1$ of an infinitely divisible subordinator $%
\left( \chi _{t};t\geq 0\right) $ with LST $\mathbf{E}\left( e^{-\lambda
\chi }\right) ^{t}$, \cite{Ber2}. In (\ref{Pcum}), the terms $\overline{\pi }%
^{-1}\left( \tau _{n}\right) $ are thus the ranked jumps of $\chi $ in its
L\'{e}vy decomposition (see \cite{Per} for example).\newline

\textbf{Reproduction step.}

Because we consider the offspring of all the $N$ initial individuals with
fitnesses $x_{n}\left( 0\right) ,$ $n=1,...,N$, the step$-1$ state of the
whole population is thus obtained from a PPP with global equivalent
occupation density 
\begin{equation}
\pi _{x_{N,\alpha }\left( 0\right) }\left( x\right) =\alpha x^{-\left(
\alpha +1\right) }\sum_{n=1}^{N}x_{n}\left( 0\right) ^{\alpha }=\alpha
x_{N,\alpha }\left( 0\right) ^{\alpha }x^{-\left( \alpha +1\right) },
\label{B2}
\end{equation}
where 
\begin{equation*}
x_{N,\alpha }\left( 0\right) :=\left( \sum_{n=1}^{N}x_{n}^{\alpha }\left(
0\right) \right) ^{1/\alpha }
\end{equation*}
is the global equivalent initial fitness of the whole population at
generation $0$ to consider$.$ We get

\begin{proposition}
In a PPP with occupation density (\ref{B1}) for the descent of each
individual with fitness $x_{n}\left( 0\right) ,$ $n=1,...,N$, the occupation
density of the population as a whole is given by (\ref{B2}), where $%
x_{N,\alpha }\left( 0\right) :=\left( \sum_{n=1}^{N}x_{n}^{\alpha }\left(
0\right) \right) ^{1/\alpha }$ is the global equivalent fitness.
\end{proposition}

\textbf{Proof:} This follows from the superposition principle of Poisson
point processes (see \cite{Kin} p. $16$). The fact that the intensity of the
superposed PPP is in the same class as the one of a single PPP descending
from $x_{n}\left( 0\right) $ is a remarkable scaling property of $\pi
_{x_{n}\left( 0\right) }\left( x\right) .$ $\diamond $\newline

In this setup therefore, $\pi _{x_{N,\alpha }\left( 0\right) }\left(
x\right) :=-\overline{\pi }_{x_{N,\alpha }\left( 0\right) }^{\prime }\left(
x\right) $ stands for the occupation density that there is a point at $x$
descending from \emph{any} of the $N$ initial individuals with fitnesses $%
x_{n}\left( 0\right) ,$ $n=1,...,N.$\newline

Note that the cumulated fitness of all first-generation offspring is 
\begin{equation*}
\sum_{n\geq 1}\overline{\pi }_{x_{N,\alpha }\left( 0\right) }^{-1}\left(
\tau _{n}\right) \overset{d}{=}x_{N,\alpha }\left( 0\right) \cdot \chi ,
\end{equation*}
where $\chi $ is given by (\ref{Pcum}).\newline

\textbf{Selection step.}

In order to model a population with fixed size over the generations, the
final state of the population at time $1$ is obtained while selecting the $N$
individuals of the whole population whose fitnesses are the largest (the
selection step), truncating therefore the latter sum to its $N$ first terms.

The whole process (including reproduction and selection steps) is then
iterated independently over the next generations.

From this definition of the process, if the $x_{n}\left( k\right) $s are the
fitnesses of the $N$ fittest individuals at generation $k$, the ordered ones 
$x_{\left( n\right) }\left( k+1\right) $ at generation $k+1$ ($x_{\left(
1\right) }>...>x_{\left( N\right) }$) descending from the whole population
at step $k$ are given by $x_{\left( n\right) }=\overline{\pi }_{x_{N,\alpha
}\left( k\right) }^{-1}\left( \tau _{n}\right) =x_{N,\alpha }\left( k\right)
\tau _{n}^{-1/\alpha }$. Here the $\tau _{n}$s are the ordered points of a
standard Poisson process on the half-line with $\tau _{1}<...<\tau
_{n}<...<\tau _{N}.$ The law of $\tau _{n}$ is thus the one of an Erlang
gamma($n$) rv with density $f_{\tau _{n}}\left( s\right)
=s^{n-1}e^{-s}/\Gamma \left( n\right) .$ Further, given $\tau _{N+1}=s$, the
probability density of $\tau _{1},...,\tau _{N}$ is 
\begin{equation*}
f_{\tau _{1},...,\tau _{N}}\left( s_{1},..,s_{N}\mid \tau _{N+1}=s\right) =%
\frac{N!}{s^{N}}1_{0<s_{1}<...<s_{N}<s}
\end{equation*}
and so 
\begin{equation*}
f_{\tau _{1},...,\tau _{N}}\left( s_{1},..,s_{N}\right)
=1_{0<s_{1}<...<s_{N}}\int_{s_{N}}^{\infty
}dse^{-s}=e^{-s_{N}}1_{0<s_{1}<...<s_{N}}.
\end{equation*}
The joint law of the ordered $x_{n}\left( k+1\right) $s is thus the one of
the images $\overline{\pi }_{x_{N,\alpha }\left( k\right) }^{-1}\left( \tau
_{n}\right) $s$,$ namely 
\begin{equation*}
f_{x_{\left( 1\right) }\left( k+1\right) ,...,x_{\left( N\right) }\left(
k+1\right) }\left( x_{1},..,x_{N}\right) =e^{-\overline{\pi }_{x_{N,\alpha
}\left( k\right) }\left( x_{N}\right) }\prod_{n=1}^{N}\pi _{x_{N,\alpha
}\left( k\right) }\left( x_{n}\right) 1_{x_{1}>...>x_{N}}.
\end{equation*}
Clearly also (with the two terms in the right-hand side term mutually
independent), 
\begin{equation*}
x_{\left( N+1\right) }\left( k+1\right) \overset{d}{=}x_{N,\alpha }\left(
k\right) x_{N+1}^{*}\left( k+1\right) ,
\end{equation*}
where $x_{N+1}^{*}\left( k+1\right) $ is the $\left( N+1\right) -$st largest
point of a PPP with occupation density $\pi \left( x\right) =\alpha
x^{-\left( \alpha +1\right) }$, $x,\alpha >0$. By the image measure theorem,
the density of $x_{N+1}^{*}\left( k+1\right) \overset{d}{=}\overline{\pi }%
^{-1}\left( \tau _{N+1}\right) $ is obtained as a power-gamma density 
\begin{equation}
f_{x_{N+1}^{*}\left( k+1\right) }\left( x\right) =\frac{\alpha }{N!}%
x^{-\left( \left( N+1\right) \alpha +1\right) }e^{-x^{-\alpha }}\text{, }x>0.
\label{B3}
\end{equation}
Next, the conditional density of each $x_{n}\left( k+1\right) $ given $%
x_{\left( N+1\right) }\left( k+1\right) =x$ is 
\begin{equation*}
f_{x_{n}\left( k+1\right) }\left( x_{n}\mid x_{\left( N+1\right) }\left(
k+1\right) =x\right) =\frac{\pi \left( x_{n}\right) }{\overline{\pi }\left(
x\right) }1_{x_{n}>x},
\end{equation*}
showing that (with the two terms in the right-hand side term mutually
independent), $x_{n}\left( k+1\right) \overset{d}{=}x_{\left( N+1\right)
}\left( k+1\right) X_{n}\left( k+1\right) $ where $X_{n}\left( k+1\right) $
is a Pareto$\left( \alpha \right) $ distributed rv with density $f\left(
x\right) =\alpha x^{-\left( \alpha +1\right) }$ on $\left( 1,\infty \right) $%
\footnote{%
We used the scaling property of Pareto$\left( \alpha \right) $ rvs $X$ on $%
\left( 1,\infty \right) $ stating that $X\mid X>a\overset{d}{=}aX.$}$.$
Putting all this together, we obtained

\begin{proposition}
Independently for each $k\geq 0$, with $x_{N+1}^{*}\left( k+1\right) $
having the power-gamma distribution (\ref{B3}) and $X_{n}\left( k+1\right) $
being Pareto$\left( \alpha \right) $ distributed and with the three
right-hand side terms being mutually independent, 
\begin{equation}
x_{n}\left( k+1\right) \overset{d}{=}x_{N,\alpha }\left( k\right)
x_{N+1}^{*}\left( k+1\right) X_{n}\left( k+1\right) \text{, }n=1,...,N,
\label{B6}
\end{equation}
indicates how to update multiplicatively the fitness of the $n-$th
individual at generation $k+1,$ when the fitnesses of the previous
generation are summarized in $x_{N,\alpha }\left( k\right) .$
\end{proposition}

We also clearly have an update of the global equivalent fitness $x_{N,\alpha
}$ from step $k$ to step $k+1$ as:

\begin{corollary}
With $x_{N,\alpha }\left( k\right) :=\left( \sum_{n=1}^{N}x_{n}\left(
k\right) ^{\alpha }\right) ^{1/\alpha }$ the global equivalent fitness of
the whole population at generation $k$, the following recursion holds 
\begin{equation}
x_{N,\alpha }\left( k+1\right) \overset{d}{=}x_{N,\alpha }\left( k\right)
x_{N+1}^{*}\left( k+1\right) \left( \sum_{n=1}^{N}X_{n}\left( k+1\right)
^{\alpha }\right) ^{1/\alpha }.  \label{B6'}
\end{equation}
\end{corollary}

In the latter sum term, each $X_{n}^{\alpha }$ is thus a Pareto$\left(
1\right) $ distributed rv with density $f\left( x\right) =x^{-2}$ on $\left(
1,\infty \right) $ and we need to sum $N$ of them independently which is
reminiscent of (\ref{P5})$.$ \newline

\textbf{Large }$N$\textbf{\ asymptotics of the }$\alpha -$\textbf{mean
fitness.}

Defining 
\begin{equation*}
\left\langle x\right\rangle _{N,\alpha }\left( k\right) :=\left( \frac{1}{N}%
\sum_{n=1}^{N}x_{n}^{\alpha }\left( k\right) \right) ^{1/\alpha },
\end{equation*}
to be the generalized (H\"{o}lder) $\alpha -$mean of the fitnesses $x_{n},$ $%
n=1,...,N,$ at generation $k,$ it follows from (\ref{B6}) that 
\begin{equation}
\left\langle x\right\rangle _{N,\alpha }\left( k+1\right) \overset{d}{=}%
\left\langle x\right\rangle _{N,\alpha }\left( k\right) x_{N+1}^{*}\left(
k+1\right) \left( \sum_{n=1}^{N}X_{n}\left( k+1\right) ^{\alpha }\right)
^{1/\alpha }.  \label{B7}
\end{equation}
By Jensen inequality, the H\"{o}lder $\alpha -$means $\left\langle
x\right\rangle _{N,\alpha }\left( k\right) $ are non-decreasing functions of 
$\alpha .$

\begin{corollary}
The $\beta -$moments of $\left\langle x\right\rangle _{N,\alpha }$ at
generation $k$ are given by
\end{corollary}

\begin{equation*}
\mathbf{E}\left( \left\langle x\right\rangle _{N,\alpha }\left( k\right)
^{\beta }\right) =\left\langle x\right\rangle _{N,\alpha }\left( 0\right)
^{\beta }\left[ \mathbf{E}\left( x_{N+1}^{*\beta }\right) \mathbf{E}\left(
\left( \sum_{n=1}^{N}X_{n}^{\alpha }\right) ^{\beta /\alpha }\right) \right]
^{k}.
\end{equation*}
\textbf{Proof:} This follows from the independence of the $x_{N+1}^{*}$s and
the $X_{n}^{\alpha }$ and their i.i.d. character within each generation $k$
and from (\ref{B7}). $\diamond $

\begin{proposition}
For large $N$, with $v_{N}:=\log \log N$%
\begin{equation*}
\frac{1}{k}\log \left\langle x\right\rangle _{N,\alpha }\left( k\right) 
\underset{k\rightarrow \infty }{\overset{a.s.}{\rightarrow }}\frac{1}{\alpha 
}v_{N}.
\end{equation*}
With 
\begin{equation*}
F_{N}\left( \beta \right) :=-\frac{\beta }{\alpha }\log \log N-\frac{\beta }{%
\alpha \log N}\left( \psi \left( 1-\beta /\alpha \right) -\log \log
N-1\right) 
\end{equation*}
and $f_{N}\left( a\right) ,$ $a<0,$ its Legendre transform, the large
deviation regime is given by 
\begin{equation*}
\frac{1}{k}\log \mathbf{P}\left( -\frac{1}{k}\log \left\langle
x\right\rangle _{N,\alpha }\left( k\right) \overset{}{\rightarrow }a\right) 
\underset{k\rightarrow \infty }{\rightarrow }f_{N}\left( a\right) \leq 0.
\end{equation*}
\end{proposition}

\textbf{Proof:} Observing that $\mathbf{E}\left( x_{N+1}^{*\beta }\right)
=\Gamma \left( N+1-\beta /\alpha \right) /\Gamma \left( N+1\right) \underset{%
N\rightarrow \infty }{\sim }N^{-\beta /\alpha }$ and applying (\ref{P5})
giving the moments of a partial sum of $N$ i.i.d. Pareto$(1)$ distributed
rvs, it follows that, with $\beta <\alpha $%
\begin{equation*}
\mathbf{E}\left( \left( \sum_{n=1}^{N}X_{n}^{\alpha }\right) ^{\beta /\alpha
}\right) \underset{N\text{ large}}{\sim }\left( N\log N\right) ^{\beta
/\alpha }\left( 1+\frac{\beta }{\alpha \log N}\left( \psi \left( 1-\beta
/\alpha \right) -\log \log N-1\right) \right) .
\end{equation*}
As a result, we get 
\begin{equation*}
\frac{\mathbf{E}\left( \left\langle x\right\rangle _{N,\alpha }\left(
k\right) ^{\beta }\right) }{\left\langle x\right\rangle _{N,\alpha }\left(
0\right) ^{\beta }}\underset{N\text{ large}}{\sim }\left( \log N\right)
^{\left( \beta k\right) /\alpha }\left( 1+\frac{\beta }{\alpha \log N}\left(
\psi \left( 1-\beta /\alpha \right) -\log \log N-1\right) \right) ^{k}.
\end{equation*}
Thus, for large $N$, with 
\begin{equation*}
F_{N}\left( \beta \right) =-\frac{\beta }{\alpha }\log \log N-\frac{\beta }{%
\alpha \log N}\left( \psi \left( 1-\beta /\alpha \right) -\log \log
N-1\right)
\end{equation*}
defining the concave thermodynamical `pressure', 
\begin{equation*}
-\frac{1}{k}\log \mathbf{E}\left( \left\langle x\right\rangle _{N,\alpha
}\left( k\right) ^{\beta }\right) \underset{k\rightarrow \infty }{\overset{}{%
\rightarrow }}F_{N}\left( \beta \right) .
\end{equation*}
Thus, with $a=F_{N}^{\prime }\left( \beta \right) <0$, by the large
deviation principle 
\begin{equation}
\frac{1}{k}\log \mathbf{P}\left( -\frac{1}{k}\log \left\langle
x\right\rangle _{N,\alpha }\left( k\right) \overset{}{\rightarrow }a\right) 
\underset{k\rightarrow \infty }{\rightarrow }f_{N}\left( a\right) \leq 0,
\label{P7}
\end{equation}
where $f_{N}\left( a\right) =\inf_{\beta <\alpha }\left( a\beta -F_{N}\left(
\beta \right) \right) $ is the concave Legendre transform of $F_{N},$ giving
the large deviation rate function of $-\log \left\langle x\right\rangle
_{N,\alpha }\left( k\right) /k$. In particular, for large $N$, with $%
v_{N}:=\log \log N$%
\begin{equation}
\frac{1}{k}\log \left\langle x\right\rangle _{N,\alpha }\left( k\right) 
\underset{k\rightarrow \infty }{\overset{a.s.}{\rightarrow }}-F_{N}^{\prime
}\left( 0\right) \sim \frac{1}{\alpha }v_{N},  \label{P7a}
\end{equation}
gives the limiting right shift of the H\"{o}lder $\alpha -$mean fitness
induced by selection effects. $\diamond $\newline

\textbf{Remarks:}

$\left( i\right) $ The limiting right-hand-side term in (\ref{P7a}),
although increasing very slowly with $N$, does not stabilize to a limit, in
contrast to other similar models \cite{BéG} of branching with selection
where, in each generation, each individual produces only two offspring with
randomly shifted fitnesses.\newline

$\left( ii\right) $ With $\alpha >0,$ let $f\left( x\right) :=x^{\alpha }>0$
define some (increasing) output map of the individuals fitnesses $x,$ with $%
x>0$. Defining 
\begin{equation*}
\left\langle f\left( x\right) \right\rangle _{N}\left( k\right) :=\frac{1}{N}%
\sum_{n=1}^{N}f\left( x_{n}\right) \left( k\right) =\frac{1}{N}%
\sum_{n=1}^{N}x_{n}\left( k\right) ^{\alpha }
\end{equation*}
to be the mean output fitness in generation $k$ of the whole population,
then, whatever $\alpha $, (\ref{P7a}) is also 
\begin{equation*}
\frac{1}{k}\log \left\langle f\left( x\right) \right\rangle _{N}\left(
k\right) \underset{k\rightarrow \infty }{\overset{a.s.}{\rightarrow }}v_{N},
\end{equation*}
interpreting the speed $v_{N}$ itself$.$ This suggests that it is of
interest to work not only on the fitnesses $x_{n}$ themselves (and their $%
\alpha -$mean $\left\langle x\right\rangle _{N,\alpha }$) but rather on some
deformed version of the fitnesses $x_{n}^{\alpha }$ (and their standard mean 
$\left\langle f\left( x\right) \right\rangle _{N}$). Clearly $\left\langle
f\left( x\right) \right\rangle _{N}$ itself obeys the recursion

\begin{equation*}
\left\langle f\left( x\right) \right\rangle _{N}\left( k+1\right) \overset{d%
}{=}\left\langle f\left( x\right) \right\rangle _{N}\left( k\right) \cdot
x_{N+1}^{*}\left( k+1\right) ^{\alpha }\cdot \sum_{n=1}^{N}X_{n}\left(
k+1\right) ^{\alpha },
\end{equation*}
deriving again from (\ref{B6}) and the definition of the equivalent global
fitness $x_{N,\alpha }\left( k\right) $.\newline

$\left( iii\right) $ Note finally that, given $x_{N,\alpha }\left( k\right) $%
, the cumulative distribution function (cdf) of the fitness $x_{\left(
1\right) }\left( k+1\right) $ of the fittest individual among the $N$
individuals at generation $k+1$ is given by 
\begin{equation*}
\mathbf{P}_{x_{N,\alpha }\left( k\right) }\left( x_{\left( 1\right) }\left(
k+1\right) \leq x\right) =\mathbf{P}\left( \overline{\pi }_{x_{N,\alpha
}\left( k\right) }^{-1}\left( \tau _{1}\right) \leq x\right) =e^{-\left(
x/x_{N,\alpha }\left( k\right) \right) ^{-\alpha }},
\end{equation*}
where $\tau _{1}$ is exp$\left( 1\right) $ distributed. Thus the fitness of
the fittest individual obeys 
\begin{equation*}
x_{\left( 1\right) }\left( k+1\right) \overset{d}{=}x_{N,\alpha }\left(
k\right) Y\left( k+1\right) ,
\end{equation*}
where $Y\left( k\right) $, $k\geq 0$ is a sequence of i.i.d. Fr\'{e}chet rvs
with cdf $\mathbf{P}\left( Y\leq x\right) =e^{-x^{-\alpha }}$. The
conditional mean given $x_{N,\alpha }\left( k\right) $ of $x_{\left(
1\right) }\left( k+1\right) $ is $x_{N,\alpha }\left( k\right) \Gamma \left(
1-1/\alpha \right) $ and its median value $x_{N,\alpha }\left( k\right)
\left( \log 2\right) ^{-1/\alpha }.$ More generally, when $N$ is large and
unconditionally, due to the recursion (\ref{B6'}) on the $x_{N,\alpha
}\left( k\right) $s$:$%
\begin{equation*}
\mathbf{E}\left( x_{\left( 1\right) }\left( k+1\right) ^{\beta }\right) \sim
x_{N,\alpha }\left( 0\right) ^{\beta }\Gamma \left( 1-\beta /\alpha \right)
e^{-kF_{N}\left( \beta \right) }.
\end{equation*}

\subsection{Genealogies}

Now we turn to the genealogies of this branching process with selection.%
\newline

\textbf{The beta}$\left( 1,1-\beta \right) $\textbf{\ and
Bolthausen-Sznitman coalescents.}

Recall that $\pi _{x_{N,\alpha }\left( k\right) }\left( x\right) :=-%
\overline{\pi }_{x_{N,\alpha }\left( k\right) }^{\prime }\left( x\right) $
is the occupation density that there would be a point (an offspring) of the
PPP at position (with fitness) $x$ at generation $k+1,$ given a global
population state $x_{N,\alpha }\left( k\right) .$

\begin{proposition}
Looking backward in time, upon scaling time using $c_{N}\sim 1/\log N$, the
genealogy of the branching model with selection is a beta$\left( 1,1-\beta
\right) $ coalescent, reducing to the Bolthausen-Sznitman coalescent if $%
\beta =0$.
\end{proposition}

\textbf{Proof:} Suppose first $\beta =0$. Given there is an offspring at $x$
at generation $k+1$, the sampling probability that it would be an offspring
of the individual with fitness $x_{n}\left( k\right) $ is thus 
\begin{equation*}
\frac{\pi _{x_{n}\left( k\right) }\left( x\right) }{\pi _{x_{N,\alpha
}\left( k\right) }\left( x\right) }=\frac{\alpha x_{n}\left( k\right)
^{\alpha }x^{-\left( \alpha +1\right) }}{\alpha x_{N,\alpha }\left( k\right)
^{\alpha }x^{-\left( \alpha +1\right) }}=\frac{x_{n}\left( k\right) ^{\alpha
}}{x_{N,\alpha }\left( k\right) ^{\alpha }},
\end{equation*}
which is independent of $x$. Observing from (\ref{B6}) and (\ref{B6'}) that 
\begin{equation*}
x_{n}\left( k\right) \overset{d}{=}x_{N,\alpha }\left( k-1\right)
x_{N+1}^{*}\left( k\right) X_{n}\left( k\right) ,\text{ }n=1,...,N,
\end{equation*}
\begin{equation*}
x_{N,\alpha }\left( k\right) \overset{d}{=}x_{N,\alpha }\left( k-1\right)
x_{N+1}^{*}\left( k\right) \left( \sum_{n=1}^{N}X_{n}\left( k\right)
^{\alpha }\right) ^{1/\alpha },
\end{equation*}
this random probability is also 
\begin{equation*}
\frac{X_{n}\left( k\right) ^{\alpha }}{\sum_{n=1}^{N}X_{n}\left( k\right)
^{\alpha }},
\end{equation*}
where, for each $k$ independently, the $X_{n}$s are i.i.d. Pareto$\left(
\alpha \right) $ distributed on $\left( 1,\infty \right) .$ Because in each
generation, parents generate offspring independently and independently of
one another, upon averaging, the probability that, at generation $k+1$, $i$
individuals share the same common ancestor is independent of $k,$ with 
\begin{equation*}
P_{i,1}^{\left( N\right) }:=\sum_{n=1}^{N}\mathbf{E}\left( \left( \frac{%
X_{n}\left( k\right) ^{\alpha }}{\sum_{n=1}^{N}X_{n}\left( k\right) ^{\alpha
}}\right) ^{i}\right) =N\mathbf{E}\left( \left( \frac{X_{1}\left( k\right)
^{\alpha }}{\sum_{n=1}^{N}X_{n}\left( k\right) ^{\alpha }}\right)
^{i}\right) .
\end{equation*}
The $X_{n}\left( k\right) $s being i.i.d. Pareto$\left( \alpha \right) $
distributed rvs, the $X_{n}\left( k\right) ^{\alpha }$s are i.i.d. Pareto$%
\left( 1\right) $ distributed rvs and we are thus back to the results of
Proposition $6$ with $\beta =0$, stating that with $c_{N}\sim 1/\log N,$ 
\begin{equation*}
c_{N}^{-1}P_{i,1}^{\left( N\right) }\underset{N\rightarrow \infty }{%
\rightarrow }\int_{0}^{1}u^{i-2}\Lambda \left( du\right) =\frac{\Gamma
\left( 2\right) \Gamma \left( i-1\right) }{\Gamma \left( i\right) }=\frac{1}{%
i-1},
\end{equation*}
where $\Lambda \overset{d}{\sim }$beta$\left( 1,1\right) ,$ uniform. We can
proceed similarly to derive the probabilities $P_{i,j}^{\left( N\right) }$
that $i$ individuals have $j<i$ parents, behaving consistently with (\ref
{P5a}) with $\beta =0$. We conclude that, whatever $\alpha $, in the large $%
N $ limit, the time-scaled genealogy of the branching model with selection
is a Bolthausen-Sznitman coalescent process, obtained while sampling from $N$
Pareto$\left( 1\right) $ i.i.d. rvs. Would the sampling probabilities $%
P_{i,j}^{\left( N\right) }$ include a $\beta -$size biasing effect on total
length $\Sigma _{N}:=\sum_{n=1}^{N}X_{n}^{\alpha }$, the genealogy of this
branching model with selection would be a full beta$\left( 1,1-\beta \right) 
$ coalescent, provided $\beta <1$. While adopting this sampling point of
view to compute the coalescence and merging probabilities, we therefore
obtain a limiting genealogical coalescent process which is independent of $%
\alpha .$ $\diamond $\newline

\textbf{Genealogies from the output PPP}.

What now if we set that the occupation density that, at generation $k+1,$
there is an offspring at $x$ descending from some individual with fitness $%
x_{n}\left( k\right) $ at generation $k$, is instead given by 
\begin{equation}
\pi _{x_{n}\left( k\right) }\left( x\right) :=x_{n}\left( k\right) x^{-2},
\label{P8a}
\end{equation}
while distorting the original occupation intensity $\pi _{x_{n}\left(
k\right) }\left( x\right) =\alpha x_{n}\left( k\right) ^{\alpha }x^{-\left(
\alpha +1\right) }$?

Then the occupation density that, at generation $k+1,$ there is an offspring
at $x$ descending from any individual of the whole population would take the
form 
\begin{equation*}
\pi _{x_{N,1}\left( k\right) }\left( x\right) :=x_{N,1}\left( k\right)
x^{-2},
\end{equation*}
where $x_{N,1}\left( k\right) :=\sum_{n=1}^{N}x_{n}\left( k\right) $ is the
cumulative fitness in generation $k$.

If this were to be the case, given there is an offspring at $x$ at time $k+1$%
, the sampling probability that it is an offspring of the individual with
fitness $x_{n}\left( k\right) $ would be, thanks to (\ref{B6}) and (\ref{B6'}%
) with $\alpha =1:$%
\begin{equation}
\frac{\pi _{x_{n}\left( k\right) }\left( x\right) }{\pi _{x_{N,1}\left(
k\right) }\left( x\right) }=\frac{x_{n}\left( k\right) x^{-2}}{x_{N,1}\left(
k\right) x^{-2}}=\frac{x_{n}\left( k\right) }{x_{N,1}\left( k\right) }%
\overset{d}{=}\frac{X_{n}\left( k\right) }{\sum_{n=1}^{N}X_{n}\left(
k\right) },  \label{P8}
\end{equation}
again independently of $x$. This probability now involves a normalized sum
of the $X_{n}$s, which are i.i.d. Pareto$\left( \alpha \right) $ distributed
and the strategy to compute the merging probabilities of the ancestral
process will be modified.

Under this hypothesis indeed, the probability that, at generation $k+1$, $i$
individuals share the same common ancestor reads 
\begin{equation*}
P_{i,1}^{\left( N\right) }:=\sum_{n=1}^{N}\mathbf{E}\left( \left( \frac{%
X_{n}\left( k\right) }{\sum_{n=1}^{N}X_{n}\left( k\right) }\right)
^{i}\right) =N\mathbf{E}\left( \left( \frac{X_{1}\left( k\right) }{%
\sum_{n=1}^{N}X_{n}\left( k\right) }\right) ^{i}\right) =:N\mathbf{E}\left(
S_{1}\left( k\right) ^{i}\right) ,
\end{equation*}
where $S_{1}$ is the normalized segment size now obtained from $N$ i.i.d.
Pareto$\left( \alpha \right) $ distributed rvs, normalized by their sum $%
\Sigma _{N}\left( k\right) :=\sum_{n=1}^{N}X_{n}\left( k\right) $. We can
proceed similarly to derive the probabilities $P_{i,j}^{\left( N\right) }$
that $i$ individuals have $j<i$ parents and we are back to the studies of
Sections $2-4$. And we can as well $\beta -$size-bias these sampling
probabilities on the total lengths $\Sigma _{N}$. We call this sampling
procedure the distorted sampling procedure.

\begin{proposition}
Looking backward in time, using a distorted size-biased sampling procedure,
the genealogy of the branching model with selection is

- a continuous-time Kingman coalescent if $\alpha \geq 2$ (upon scaling time
with $c_{N}\propto 1/N$ if $\alpha >2$ or $c_{N}\propto \log N/N$ if $\alpha
=2$).

- a continuous-time beta$\left( 2-\alpha ,\alpha -\beta \right) $ coalescent
if $\alpha \in \left( 1,2\right) $ and $\beta <\alpha $ (upon scaling time
with $c_{N}\propto N^{-\left( \alpha -1\right) }$).

- a continuous-time beta$\left( 1,1-\beta \right) $ coalescent if $\alpha =1$
and $\beta <1$ (upon scaling time with $c_{N}\propto 1/\log N$).

- a discrete-time Poisson-Dirichlet$\left( \alpha ,-\beta \right) $
coalescent if $\alpha \in \left( 0,1\right) $ and $\beta <\alpha $.
\end{proposition}

\textbf{Proof:} It remains to interpret the distorted size-biased sampling
procedure which is proposed to compute the merging probabilities of the
ancestral process: Assume that the fitness dependent PPP describing the
descent of an individual with fitness $x_{n}\left( k\right) $ is now the
output image of the original one, given by the canonical application $%
f_{x_{n}\left( k\right) }\left( x\right) =x_{n}\left( k\right) \left(
x/x_{n}\left( k\right) \right) ^{\alpha },$ $f_{x_{n}\left( k\right) }:\Bbb{R%
}_{+}\rightarrow \Bbb{R}_{+}.$ Its gives rise to a new PPP with distorted
intensity $\pi _{f_{x_{n}\left( k\right) }}\left( dx\right) =x_{n}\left(
k\right) x^{-2}dx$, the image measure of $\pi _{x_{n}\left( k\right) }\left(
dx\right) =\pi _{x_{n}\left( k\right) }\left( x\right) dx=\alpha x_{n}\left(
k\right) ^{\alpha }x^{-\left( \alpha +1\right) }dx$ by $f_{x_{n}\left(
k\right) }$ (as a result of the classical Campbell formula for PPPs, see 
\cite{Kin} p. $28$). This was the starting point in (\ref{P8a}). The skewed
computation in (\ref{P8}) of the probability that some offspring is
descending from one individual with fitness $x_{n}\left( k\right) $ is thus
based not on the original PPP attached to $x_{n}\left( k\right) $ but rather
on a deformed version of it through $f_{x_{n}\left( k\right) }.$ We note
that $f_{x_{n}\left( k\right) }\left( x\right) $, as a function of the two
arguments $\left( x_{n}\left( k\right) ,x\right) $ is homogeneous with $%
f_{\lambda ^{a}x_{n}\left( k\right) }\left( \lambda ^{b}x\right) =\lambda
^{a\left( 1-\alpha \right) +b\alpha }f_{x_{n}\left( k\right) }\left(
x\right) $, $\lambda >0,$ leading obviously, if $\lambda =x_{n}\left(
k\right) ^{-1/a}$, $b/a=1$ and $f\left( x\right) :=f_{1}\left( x\right)
=x^{\alpha }$ to: $f_{x_{n}\left( k\right) }\left( x\right) =x_{n}\left(
k\right) f\left( x/x_{n}\left( k\right) \right) $. The function $f\left(
x\right) =x^{\alpha }$ is the output fitness function introduced in
Subsection $5.1,$ Remark$\left( ii\right) $.

Using this distorted size-biased sampling procedure therefore, following the
introductory arguments, the full class of the Pareto-coalescents (described
in Sections $2-4$) are obtained.

Suppose for instance $\alpha \in \left( 1,2\right) $, $\beta =0$. Based on
the previous computations of Sections $2-4$, we conclude that the large $N$
distorted genealogy of the branching model with selection coincides (upon
scaling time correspondingly: $k\rightarrow \left[ t/c_{N}\right] $) with a
beta$\left( 2-\alpha ,\alpha \right) $ coalescent process, obtained while
sampling from $N$ Pareto$\left( \alpha \right) $ i.i.d. rvs, normalized by
their sum. Would this probability involve a $\beta -$size biasing effect,
the genealogy of this branching model with selection is identified to a beta$%
\left( 2-\alpha ,\alpha -\beta \right) $ coalescent, $\beta <\alpha $.

If $\alpha \in \left[ 0,1\right) ,$ $\beta <\alpha $, the obtained large $N$
genealogical coalescent will coincide with the discrete-time-$k$
Poisson-Dirichlet coalescent with parameters $\alpha $ and $-\beta .$ Only
when $\alpha =1$ do we get as in \cite{BD} (upon scaling time
logarithmically with $N$) the Bolthausen-Sznitman coalescent ($\beta =0$) or
more generally the beta$\left( 1,1-\beta \right) $ coalescent, provided $%
\beta <1$. $\diamond $\newline

\textbf{Acknowledgments:} The author acknowledges partial support from the
ANR Mod\'{e}lisation Al\'{e}atoire en \'{E}cologie, G\'{e}n\'{e}tique et 
\'{E}volution (ANR-Man\`{e}ge- 09-BLAN-0215 project) and from the labex
MME-DII (Mod\`{e}les Math\'{e}matiques et \'{E}conomiques de la Dynamique,
de l' Incertitude et des Interactions). The author is also indebted to his
referees for pointing out some errors in an earlier version of the draft and
for encouraging him to write down a more concise and complete version.

\end{document}